# Partial Gâteaux and Fréchet Derivatives and Applications to Variational Analysis

Jinlu Li

Department of Mathematics
Shawnee State University
Portsmouth, Ohio 45662 USA
jli@shawnee.edu

**Abstract.** In this paper, we define partial Gâteaux and Fréchet derivatives for multivariable and single-valued mappings between Banach spaces. We will prove some properties of partial Gâteaux and Fréchet derivatives. By these definitions, we find the explicit formulas for some polynomial type operators with two variables from $l_p \times l_p$ to $l_p$ with respect to each variable. Then, we will introduce the concepts of generalized partially critical points and partial ordered extrema in partially ordered Banach spaces. By these concepts, we will investigate the connection between generalized partially critical points and partial ordered extrema of two-variable and single-valued mappings in partially ordered Banach spaces. These results extend the connection between critical points and extrema of real valued functions in calculus. We will give some applications of partial Gâteaux and Fréchet derivatives to ordered optimizations and variational inequality problems between partially ordered Banach spaces. Finally, we study the connection between partial Gâteaux and Fréchet derivatives and ordered monotone of two-variable and single-valued mappings in partially ordered Banach spaces.

## 1. Introduction and Preliminary

The theory of differentiation in Banach spaces has been developed rapidly, which plays a prominent role in modern analysis. In contrast with differentiation in ordinary calculus, the most popular and most useful concepts of differentiation in Banach spaces are Gâteaux (directional) differentiability and Fréchet differentiability, which are reviewed below.

Let $(X, \|\cdot\|_X)$ and $(Z, \|\cdot\|_Z)$ be Banach spaces with null elements $\theta_X$ and $\theta_Z$, respectively. Let $D$ be a nonempty convex and open subset of $X$ and let $\bar{x} \in D$. Let $T: D \to Z$ be a single-valued mapping.

(Gâteaux directional differentiability of $T$). Let $v \in X$ with $v \neq \theta_X$. If there is a vector in $Z$, which is denoted by $T'(\bar{x}, v)$ such that

$$T'(\bar{x}, v) = \lim_{t \to 0} \frac{T(\bar{x} + tv) - T(\bar{x})}{t}, \tag{1.1}$$

then, $T$ is said to be Gâteaux directionally differentiable at point $\bar{x}$ along direction $v$. The point $T'(\bar{x}, v) \in Z$ is called the Gâteaux directional derivative of $T$ at point $\bar{x}$ along direction $v$. Furthermore, if $T$ is Gâteaux directionally differentiable at point $\bar{x}$ along every direction $v \in X \backslash \{\theta_X\}$, then $T$ is said to be Gâteaux differentiable at $\bar{x}$ and the Gâteaux derivative of $T$ at point $\bar{x}$ is denoted by

$$T'(\bar{x})(v) = \lim_{t \to 0} \frac{T(\bar{x}+tv)-T(\bar{x})}{t}, \text{ for any } v \in X\backslash\{\theta_X\}. \quad (1.2)$$

(Fréchet differentiability of $T$). If there is a continuous and linear mapping $\nabla T(\bar{x}): X \to Z$ such that,

$$\lim_{\substack{X \\ u \to \theta_X}} \frac{T(\bar{x}+u)-T(\bar{x})-\nabla T(\bar{x})(u)}{\|u\|_X} = \theta_Z, \quad (1.3)$$

then $T$ is said to be Fréchet differentiable at $\bar{x}$ and $\nabla T(\bar{x})$ is called the Fréchet derivative of $T$ at $\bar{x}$. It is well-known that in Banach spaces (or in normed vector spaces), for single-valued mappings, Fréchet differentiability is stronger than Gâteaux differentiability. That is,

$$T \text{ is Fréchet differentiable at } \bar{x} \quad \Longrightarrow \quad T \text{ is Gâteaux differentiable at } \bar{x}. \quad (1.4)$$

More precisely speaking, the connection between Gâteaux differentiability and Fréchet differentiability in Banach spaces is that if $T$ is Fréchet differentiable at $\bar{x}$, then $T$ is also Gâteaux differentiable at $\bar{x}$ and

$$T'(\bar{x})(v) = \nabla T(\bar{x})(v), \text{ for any } v \in X\backslash\{\theta_X\}. \quad (1.5)$$

However, in [15] two counter examples are provided to show that the inverse of (1.4) does not hold:

$$T \text{ is Gâteaux differentiable at } \bar{x} \quad \not\Rightarrow \quad T \text{ is Fréchet differentiable at } \bar{x}. \quad (1.6)$$

Both Gâteaux (directional) differentiation and Fréchet differentiation in Banach spaces have been widely applied to many fields in both pure and applied mathematics (See [1−4, 6−10, 20]), in particular in game theory (See [18−20]), optimization theory (See [11−12, 18−19]), and so forth.

It is well-known that, in real analysis, partial derivatives of real valued functions appear in many calculus-based problems with more than one variable in Euclidean spaces, such as: partial differential equations, optimization problems; economics (see [5]), and so forth.

This motivates us to study partial Gâteaux and Fréchet derivatives of single-valued mappings with more than one variable in Banach spaces. To extend the concepts of Gâteaux and Fréchet derivatives of single-valued mappings with one variable, in section 2, we will define partial Gâteaux and Fréchet derivatives of single-valued and multivariable mappings in Banach spaces. Some properties of partial Gâteaux and Fréchet derivatives will be proved. In [13−15], the exact solutions of Gâteaux and Fréchet derivatives of some single-valued mappings with one variable have been calculated. In contrast with single-valued mappings, in section 3, by these definitions of partial Gâteaux and Fréchet derivatives, we find the explicit formulas for the partial Gâteaux and Fréchet derivatives of some polynomial type operators from $l_p \times l_p$ to $l_p$ with respect to each variable. In section 4, by partial Gâteaux and Fréchet derivatives, we define generalized partially critical points of multivariable mappings with respect to a given variable. We will generalize properties (1.4) and (1.6) of single-valued mappings to multivariable mappings.

In sections 4, 5, and 6, we will give some applications of partial Gâteaux and Fréchet derivatives of multivariable mappings to partially ordered extrema and partially ordered optimization problems with respect to some power mappings, duality mappings.

In sections 7 and 8, we study some applications of partial Gâteaux and Fréchet derivatives to ordered variational inequality problems and ordered monotone properties of multivariable mappings between partially ordered Banach spaces.

## 2. Partial Gâteaux and Fréchet Differentiation in Banach Spaces

## 2.1 Definitions and Properties of Partial Gâteaux and Fréchet Derivatives in Banach Spaces

In this section, we extend the concepts of Gâteaux (directional) differentiability defined by (1.1) and (1.2) and Fréchet differentiability defined by (1.3) of single variable and single-valued mappings to partial Gâteaux and Fréchet differentiability of multivariable and single-valued mappings between Banach spaces, which will be used to study ordered optimization problems and ordered variational inequality problems in partially ordered Banach spaces. In this paper, conventionally, for clarity and simplicity of notation, we only consider two-variable and single-valued mappings, which can be similarly extended to arbitrary multivariable mappings.

Let $(X, \|\cdot\|_X)$, $(Y, \|\cdot\|_Y)$ and $(Z, \|\cdot\|_Z)$ be Banach spaces with null elements $\theta_X$, $\theta_Y$ and $\theta_Z$ in $X$, $Y$ and $Z$, respectively. Let $\theta_{X,Z}: X \to Z$ be the zero continuous and linear mapping from $X$ to $Z$, which satisfies that

$$\theta_{X,Z}(v) = \theta_Z, \text{for any } v \in X.$$

Similar to $\theta_{X,Z}$, let $\theta_{Y,Z}: X \to Z$ be the zero continuous and linear mapping from $Y$ to $Z$ satisfying

$$\theta_{Y,Z}(u) = \theta_Z, \text{for any } u \in Y.$$

Let $A$ and $B$ be nonempty convex and open subsets in $X$ and $Y$, respectively. Let $F: A\times B \to Z$ be a single-valued mapping. Let $(\bar{x}, \bar{y}) \in A \times B$. We first define partial Gâteaux (directional) differentiability of $F$ at point $(\bar{x}, \bar{y})$.

**Definition 2.1**. Let $F: A\times B \to Z$ be a single-valued mapping. Let $(\bar{x}, \bar{y}) \in A \times B$. Given $v \in X\backslash\{\theta_X\}$. If there is $\frac{\partial F(x,y)}{\partial x}(\bar{x}, \bar{y}, v) \in Z$ such that

$$\frac{\partial F(x,y)}{\partial x}(\bar{x}, \bar{y}, v) = \lim_{t\to 0}\frac{F(\bar{x}+tv,\bar{y})-F(\bar{x},\bar{y})}{t}, \tag{2.1}$$

then, $F$ is said to be partially Gâteaux directionally differentiable with respect to the first variable at point $(\bar{x}, \bar{y})$ along the direction $v \in X\backslash\{\theta_X\}$ and $\frac{\partial F(x,y)}{\partial x}(\bar{x}, \bar{y})(v)$ is called the partially Gâteaux directional derivative. If (2.1) holds for each $v \in X\backslash\{\theta_X\}$, then $F$ is said to be partially Gâteaux differentiable with respect to the first variable at point $(\bar{x}, \bar{y})$ with partial Gâteaux derivative $\frac{\partial F(x,y)}{\partial x}(\bar{x}, \bar{y})$ and

$$\frac{\partial F(x,y)}{\partial x}(\bar{x}, \bar{y})(v) = \lim_{t\to 0}\frac{F(\bar{x}+tv,\bar{y})-F(\bar{x},\bar{y})}{t}, \text{ for any } v \in X\backslash\{\theta_X\}. \tag{2.2}$$

Given $w \in Y\backslash\{\theta_Y\}$. We can similarly define the partially Gâteaux directional derivative of $F$ with respect to the second variable at point $(\bar{x}, \bar{y})$ along the direction $w$ as follows

$$\frac{\partial F(x,y)}{\partial y}(\bar{x}, \bar{y}, w) = \lim_{t\to 0}\frac{F(\bar{x},\bar{y}+tw)-F(\bar{x},\bar{y})}{t}. \tag{2.3}$$

If (2.3) holds for any $w \in Y\backslash\{\theta_Y\}$, then $F$ is said to be partially Gâteaux differentiable with respect to the second variable at point $(\bar{x}, \bar{y})$ with partial Gâteaux derivative $\frac{\partial F(x,y)}{\partial y}(\bar{x}, \bar{y})(w)$ and

$$\frac{\partial F(x,y)}{\partial y}(\bar{x}, \bar{y})(w) = \lim_{t\to 0}\frac{F(\bar{x},\bar{y}+tw)-F(\bar{x},\bar{y})}{t}, \text{ for any } w \in Y\backslash\{\theta_Y\}. \tag{2.4}$$

Partial Gâteaux (directional) derivatives have the following linearity properties, which is similar to the linearity properties of Gâteaux (directional) derivatives.

**Theorem 2.2**. *Let $F, G$: $A\times B \to Z$ be single-valued mappings. Let $(\bar{x}, \bar{y}) \in A \times B$ and $v \in X\backslash\{\theta_X\}$. Let $a$ and b be scalars. Suppose that $F$ and $G$ are partially Gâteaux directionally differentiable with respect to the first variable at $(\bar{x}, \bar{y})$ along direction v. Then $aF + bG$ is also partially Gâteaux directionally differentiable with respect to the first variable at $(\bar{x}, \bar{y})$ along direction v and*

$$\frac{\partial(aF+bG)(x,y)}{\partial x}(\bar{x}, \bar{y}, v) = a\frac{\partial F(x,y)}{\partial x}(\bar{x}, \bar{y}, v) + b\frac{\partial G(x,y)}{\partial x}(\bar{x}, \bar{y}, v). \tag{2.5}$$

*Furthermore, if both $F$ and $G$ are partially Gâteaux differentiable with respect to the first variable at $(\bar{x}, \bar{y})$ along direction v, then so is $aF + bG$ and*

$$\frac{\partial(aF+bG)(x,y)}{\partial x}(\bar{x}, \bar{y})(v) = a\frac{\partial F(x,y)}{\partial x}(\bar{x}, \bar{y})(v) + b\frac{\partial G(x,y)}{\partial x}(\bar{x}, \bar{y})(v), \text{ for any } v \in X\backslash\{\theta_X\}. \tag{2.6}$$

*The partial Gâteaux* (*directional*) *derivatives with respect to the second variable of $F$ and $G$ have the properties same with* (2.5) *and* (2.6).

*Proof*. The proof of this theorem is straight forward and it is omitted here. □

We write notations for partial Gâteaux (directional) derivative of $F$ with respect to the first variable:

$$\frac{\partial F(x,y)}{\partial x}(\bar{x}, \bar{y}, v) = \partial_1 F(\bar{x}, \bar{y}, v) \quad \text{and} \quad \frac{\partial F(x,y)}{\partial x}(\bar{x}, \bar{y})(v) = \partial_1 F(\bar{x}, \bar{y})(v).$$

Similar to the first variable, with respect to the second variable, we write

$$\frac{\partial F(x,y)}{\partial y}(\bar{x}, \bar{y}, u) = \partial_2 F(\bar{x}, \bar{y}, u) \quad \text{and} \quad \frac{\partial F(x,y)}{\partial y}(\bar{x}, \bar{y})(u) = \partial_2 F(\bar{x}, \bar{y})(u).$$

Now, we define partial Fréchet differentiability of $F$ at point $(\bar{x}, \bar{y})$.

**Definition 2.3**. Let $F$: $A\times B \to Z$ be a single-valued mapping. Let $(\bar{x}, \bar{y}) \in A \times B$. If there is a continuous and linear mapping $\nabla_1 F(\bar{x}, \bar{y})$: $X \to Z$ such that,

$$\lim_{u \overset{X}{\to} \theta_X} \frac{F(\bar{x}+u, \bar{y}) - F(\bar{x}, \bar{y}) - \nabla_1 F(\bar{x}, \bar{y})(u)}{\|u\|_X} = \theta_Z, \tag{2.7}$$

then $F$ is said to be partially Fréchet differentiable with respect to the first variable at point $(\bar{x}, \bar{y})$ and $\nabla_1 F(\bar{x}, \bar{y})$ is called the partial Fréchet derivative. We can similarly define the partial Fréchet derivative $\nabla_2 F(\bar{x}, \bar{y})$: $Y \to Z$ of $F$ with respect to the second variable at point $(\bar{x}, \bar{y})$ as follows

$$\lim_{w \overset{Y}{\to} \theta_Y} \frac{F(\bar{x}, \bar{y}+w) - F(\bar{x}, \bar{y}) - \nabla_2 F(\bar{x}, \bar{y})(w)}{\|w\|_Y} = \theta_Z. \tag{2.8}$$

Similarly to Theorem 2.2, we can prove the linearity of partial Fréchet derivatives.

**Theorem 2.4**. *Let $F, G$: $A\times B \to Z$ be single-valued mappings. Let $(\bar{x}, \bar{y}) \in A \times B$ and let $a$ and b be scalars. Suppose that $F$* and *$G$ are partially Fréchet differentiable with respect to the first variable at $(\bar{x}, \bar{y})$. Then $aF + bG$ is also partially Fréchet differentiable with respect to the first variable at $(\bar{x}, \bar{y})$ and*

$$\nabla_1(aF + bG)(\bar{x}, \bar{y}) = a\nabla_1 F(\bar{x}, \bar{y}) + b\nabla_1 G(\bar{x}, \bar{y}). \tag{2.9}$$

*The partial Fréchet derivatives with respect to the second variable of $F$ and $G$ have the same properties as* (2.9).

*Proof*. The proof of this theorem is straight forward and it is omitted here. □

Similar to property (1.4) that Fréchet differentiability of single-valued mappings is stronger than Gâteaux differentiability, for multivariable single-valued mappings between Banach spaces, partial Fréchet differentiability is also stronger than partial Gâteaux differentiability. This property is proved by the following theorem.

**Theorem 2.5**. *Let $F: A\times B \to Z$ be a single-valued mapping. Let $(\bar{x}, \bar{y}) \in A \times B$. If F is partially Fréchet differentiable with respect to the first* (*second*) *variable at point $(\bar{x}, \bar{y})$ with the partial Fréchet derivative* $\nabla_1 F(\bar{x}, \bar{y})$ ($\nabla_2 F(\bar{x}, \bar{y})$), *then F is partially Gâteaux differentiable with respect to the first* (*second*) *variable at point $(\bar{x}, \bar{y})$ and*

$$\partial_1 F(\bar{x}, \bar{y})(v) = \nabla_1 F(\bar{x}, \bar{y})(v), \text{ for any } v \in X\backslash\{\theta_X\},$$

*and*
$$\partial_2 F(\bar{x}, \bar{y})(u) = \nabla_2 F(\bar{x}, \bar{y})(u), \text{ for any } u \in Y\backslash\{\theta_Y\}. \tag{2.10}$$

*Proof*. Suppose that $F$ is partially Fréchet differentiable with respect to the first variable at point $(\bar{x}, \bar{y})$ with the partial Fréchet derivative $\nabla_1 F(\bar{x}, \bar{y})$. For any $v \in X\backslash\{\theta_X\}$, by (2.7), we have

$$\lim_{t\to 0}\left(\frac{F(\bar{x}+tv,\bar{y})-F(\bar{x},\bar{y})}{t} - \nabla_1 F(\bar{x}, \bar{y})(v)\right)$$

$$= \lim_{t\to 0}\frac{F(\bar{x}+tv,\bar{y})-F(\bar{x},\bar{y})-\nabla_1 F(\bar{x},\bar{y})(tv)}{t}$$

$$= \|v\|_X \lim_{t\to 0}\frac{F(\bar{x}+tv,\bar{y})-F(\bar{x},\bar{y})-\nabla_1 F(\bar{x},\bar{y})(tv)}{t\|v\|_X}$$

$$= \|v\|_X \lim_{\|tu\|_X\to 0}\frac{F(\bar{x}+tv,\bar{y})-F(\bar{x},\bar{y})-\nabla_1 F(\bar{x},\bar{y})(tv)}{\|tv\|_X}$$

$$= \theta_Z.$$

This implies that $F$ is partially Gâteaux differentiable with respect to the first variable at point $(\bar{x}, \bar{y})$ and

$$\partial_1 F(\bar{x}, \bar{y})(v) = \nabla_1 F(\bar{x}, \bar{y})(v), \text{ for any } v \in X\backslash\{\theta_X\}.$$

We can similarly prove

$$\partial_2 F(\bar{x}, \bar{y})(w) = \nabla_2 F(\bar{x}, \bar{y})(w), \text{ for any } w \in Y\backslash\{\theta_Y\}.$$ □

By using partial Gâteaux derivatives, we define generalized partially critical points of multivariable and single-valued mappings between Banach spaces.

**Definition 2.6**. Let $F: A\times B \to Z$ be a single-valued mapping. Let $(\bar{x}, \bar{y}) \in A \times B$. If $F$ is partially Gâteaux differentiable with respect to the first variable at point $(\bar{x}, \bar{y})$ and

$$\partial_1 F(\bar{x}, \bar{y})(v) = \theta_Z, \text{ for any } v \in X\backslash\{\theta_X\}, \tag{2.11}$$

then, this point $(\bar{x}, \bar{y})$ is called a generalized partially critical point of $F$ with respect to the first variable. This point $(\bar{x}, \bar{y})$ is a generalized partially critical point of $F$ with respect to the second variable can be similarly defined by

$$\partial_2 F(\bar{x}, \bar{y})(u) = \theta_Z, \text{ for any } u \in Y\backslash\{\theta_Y\}. \tag{2.12}$$

More strictly, if $(\bar{x}, \bar{y})$ is simultaneously generalized partially critical points of $F$ with respect to each variable, which satisfies that,

$$\partial_1 F(\bar{x}, \bar{y})(v) = \theta_Z, \text{ for any } v \in X \backslash \{\theta_X\},$$

and
$$\partial_2 F(\bar{x}, \bar{y})(u) = \theta_Z, \text{ for any } u \in Y \backslash \{\theta_Y\}, \tag{2.13}$$

then $(\bar{x}, \bar{y})$ is called a generalized critical point of $F$.

By the connection between partial Gâteaux and partial Fréchet differentiability given in Theorem 2.5, partial Fréchet derivatives also define the generalized partially critical points.

**Lemma 2.7**. *Let $F: A \times B \to Z$ be a single-valued mapping. Let $(\bar{x}, \bar{y}) \in A \times B$. If F is partially Fréchet differentiable with respect to the first (second) variable at point $(\bar{x}, \bar{y})$ such that*

$$\nabla_1 F(\bar{x}, \bar{y}) = \theta_{X,Z} \; (\nabla_2 F(\bar{x}, \bar{y}) = \theta_{Y,Z}).$$

*then, $(\bar{x}, \bar{y})$ is a generalized partially critical point of F with respect to the first (second) variable.*

*Proof*. By the conditions of this lemma and Theorem 2.5, $F$ is partially Gâteaux differentiable with respect to the first (second) variable at point $(\bar{x}, \bar{y})$. Then, by (2.10) in Theorem 2.5, we have

$$\partial_1 F(\bar{x}, \bar{y})(v) = \nabla_1 F(\bar{x}, \bar{y})(v) = \theta_{X,Z}(v) = \theta_Z, \text{ for any } v \in X \backslash \{\theta_X\},$$

and
$$\partial_2 F(\bar{x}, \bar{y})(u) = \nabla_2 F(\bar{x}, \bar{y})(u) = \theta_{Y,Z}(u) = \theta_Z, \text{ for any } u \in Y \backslash \{\theta_Y\}.$$ □

We define the following notations to name the sets of generalized partially critical points

**Notations 2.8**. Let $F: A \times B \to Z$ be a single-valued mapping. We define the following sets.

(i) The set of generalized partially critical points of $F$ with respect to the first variable:

$$\partial_1^C F := \{(\bar{x}, \bar{y}) \in A \times B: \partial_1 F(\bar{x}, \bar{y})(v) = \theta_Z, \text{for any } v \in X \backslash \{\theta_X\}\};$$

(ii) The set of generalized partially critical points of $F$ with respect to the second variable:

$$\partial_2^C F := \{(\bar{x}, \bar{y}) \in A \times B: \partial_2 F(\bar{x}, \bar{y})(u) = \theta_Z, \text{for any } u \in Y \backslash \{\theta_Y\}\};$$

(iii) The set of generalized critical points of $F$:

$$\partial^C F := \{(\bar{x}, \bar{y}) \in A \times B: \partial_1 F(\bar{x}, \bar{y})(v) = \theta_Z \text{ and } \partial_2 F(\bar{x}, \bar{y})(u) = \theta_Z, \text{for any } v \in X \backslash \{\theta_X\}, u \in Y \backslash \{\theta_Y\}\}.$$

(iv) The set of generalized partially critical points of $F$ with respect to partial Fréchet derivatives:

$$\nabla_1^C F := \{(\bar{x}, \bar{y}) \in A \times B: \nabla_1 F(\bar{x}, \bar{y}) = \theta_{X,Z}\};$$

$$\nabla_2^C F := \{(\bar{x}, \bar{y}) \in A \times B: \nabla_2 F(\bar{x}, \bar{y}) = \theta_{Y,Z}\};$$

$$\nabla^C F := \{(\bar{x}, \bar{y}) \in A \times B: \nabla_1 F(\bar{x}, \bar{y}) = \theta_{X,Z} \text{ and } \nabla_2 F(\bar{x}, \bar{y}) = \theta_{Y,Z}\}.$$

**Lemma 2.9**. Let $F: A \times B \to Z$ be a single-valued mapping. We have

$$\nabla_1^C F \subseteq \partial_1^C F, \quad \nabla_2^C F \subseteq \partial_2^C F \quad \text{and} \quad \nabla^C F \subseteq \partial^C F.$$

*Proof*. By Lemma 2.7, the proof of this lemma is straight forward and it is omitted here. □

Next, we consider some special two-variable single-valued mappings.

**Lemma 2.10**. *Let $f: A \to Z$ and $g: B \to Z$ be single-valued mappings. Let $F: A \times B \to Z$ be a single-valued mapping defined by*

$$F(x, y) = f(x) + g(y), \text{ for any } (x, y) \in A \times B.$$

(i) *$F$ is partially Gâteaux differentiable with respect to the first variable at point $(\bar{x}, \bar{y})$ if and only if $f$ is Gâteaux differentiable at point $\bar{x}$. In this case,*

$$\partial_1 F(\bar{x}, \bar{y})(v) = f'(\bar{x})(v), \text{ for any } v \in X \backslash \{\theta_X\}.$$

(ii) *$F$ is partially Gâteaux differentiable with respect to the second variable at point $(\bar{x}, \bar{y})$ if and only if $g$ is Gâteaux differentiable at point $\bar{y}$. In this case,*

$$\partial_2 F(\bar{x}, \bar{y})(u) = g'(\bar{y})(u), \text{ for any } u \in Y \backslash \{\theta_Y\}.$$

(iii) *$F$ is partially Fréchet differentiable with respect to the first variable at point $(\bar{x}, \bar{y})$ if and only if $f$ is Fréchet differentiable at point $\bar{x}$. In this case,*

$$\nabla_1 F(\bar{x}, \bar{y})(v) = \nabla f(\bar{x})(v), \text{ for any } v \in X \backslash \{\theta_X\}.$$

(iv) *$F$ is partially Fréchet differentiable with respect to the second variable at point $(\bar{x}, \bar{y})$ if and only if $g$ is Fréchet differentiable at point $\bar{y}$. In this case,*

$$\nabla_2 F(\bar{x}, \bar{y})(v) = \nabla g(\bar{y})(u), \text{ for any } u \in Y \backslash \{\theta_Y\}.$$

*Proof*. Here, we write

$$F(x, y) = f(x) + g(y), \text{ for any } (x, y) \in A \times B.$$

Then, by definitions, the proof of this lemma is straight forward and it is omitted here. □

### 2.2. Counter Examples

Theorem 2.5 shows that for single-valued mappings, partial Fréchet differentiability is stronger than partial Gâteaux differentiability. That is, with respect to a given variable,

$$F \text{ is partially Fréchet differentiable at } (\bar{x}, \bar{y}) \implies F \text{ is partially Gâteaux differentiable at } (\bar{x}, \bar{y}). \quad (2.14)$$

Similarly to (1.6), we provide two counter examples below to show that the inverse of (2.14) does not hold. That is, with respect to a given variable,

$$F \text{ is partially Gâteaux differentiable at } (\bar{x}, \bar{y}) \nRightarrow F \text{ is partially Fréchet differentiable at } (\bar{x}, \bar{y}). \quad (2.15)$$

More precisely, in this subsection, at first, we will use the counter example given in Theorem 3.1 in [18] to construct some counter examples to demonstrate the property in (2.15). We first review some notations used in Theorem 3.1 in [18].

Let $(X, \|\cdot\|_X)$ be a uniformly convex and uniformly smooth Banach space with topological dual space $(X^*, \|\cdot\|_{X^*})$. Let $\langle\cdot, \cdot\rangle_X$ be the pairing between $X^*$ and $X$. Let $J: X \to X^*$ be the normalized duality mapping on $X$. Let $\mathbb{B}$ and $\mathbb{S}$ denote the closed unit ball and unit sphere in $X$, respectively. Let $P_{\mathbb{B}}: X \to \mathbb{B}$ be the standard metric projection. By Lemma 2.8 in [18], $P_{\mathbb{B}}$ has the following representation:

$$P_{\mathbb{B}}(x) = \begin{cases} x, & \text{for any } x \in \mathbb{B}, \\ \frac{1}{\|x\|_X}x, & \text{for any } x \in X\backslash\mathbb{B}. \end{cases}$$

In [17−18], we introduced the following notations, which will be used for the construction of the desired counter example. For any $x \in \mathbb{S}$, we write

(a) $x^{\uparrow} = \{v \in X\backslash\{\theta\}$: there is $\delta > 0$ such that $\|x + tv\|_X > 1$, for all $t \in (0, \delta)\}$;
(b) $x^{\downarrow} = \{v \in X\backslash\{\theta\}$: there is $\delta > 0$ such that $\|x + tv\|_X < 1$, for all $t \in (0, \delta)\}$.

In part (iii) of Theorem 3.1 in [18], let $r = 1$, then we have the following special case.

(Partial of Theorem 3.1 in [18]). *Let X be a uniformly convex and uniformly smooth Banach space. The metric projection $P_{\mathbb{B}}: X \to \mathbb{B}$ has the following properties.*

(I) *$P_{\mathbb{B}}$ is not Fréchet differentiable at any $\bar{x} \in \mathbb{S}$; $\nabla P_{\mathbb{B}}(\bar{x})$ does not exist, for any $\bar{x} \in \mathbb{S}$.*

(II) *$P_{\mathbb{B}}$ is Gâteaux differentiable on $\mathbb{S}$ satisfying that, for every point $\bar{x} \in \mathbb{S}$, we have*

(a) $P'_{\mathbb{B}}(\bar{x})(w) = w - \langle J(\bar{x}), w\rangle_X \bar{x}$, *if* $w \in \bar{x}^{\uparrow}$;
(b) $P'_{\mathbb{B}}(\bar{x})(\bar{x}) = \theta_X$ (this is the special case in (a) with $w = \bar{x}$);
(c) $P'_{\mathbb{B}}(\bar{x})(w) = w$, *if* $w \in \bar{x}^{\downarrow}$.

**Theorem 2.11**. *Let $(X, \|\cdot\|_X)$ be a uniformly convex and uniformly smooth Banach space. Define $F: X\times X \to X$ as follows*:

$$F(x, y) = P_{\mathbb{B}}(x) + P_{\mathbb{B}}(y), \text{ for any } (x, y) \in X\times X.$$

*Let $(\bar{x}, \bar{y}) \in \mathbb{S} \times \mathbb{S}$ be arbitrarily given. Then*

(I) *F is not partially Fréchet differentiable with respect to any variable at $(\bar{x}, \bar{y})$. That is,*

*neither $\nabla_1 F(\bar{x}, \bar{y})$, nor $\nabla_2 F(\bar{x}, \bar{y})$ exists.*

(II) *F is partially Gâteaux differentiable with respect to each variable at $(\bar{x}, \bar{y})$ and*

(i) *With respect to the first variable, we have*

(a) $\partial_1 F(\bar{x}, \bar{y})(v) = P'_{\mathbb{B}}(\bar{x})(v) = v - \langle J(\bar{x}), v\rangle_X \bar{x}$, *if* $v \in \bar{x}^{\uparrow}$;
(b) $\partial_1 F(\bar{x}, \bar{y})(\bar{x}) = P'_{\mathbb{B}}(\bar{x})(\bar{x}) = \theta_X$ *(this is the special case in* (a) *with* $v = \bar{x}$);
(c) $\partial_1 F(\bar{x}, \bar{y})(v) = P'_{\mathbb{B}}(\bar{x})(v) = w$, *if* $v \in \bar{x}^{\downarrow}$.

(ii) *With respect to the second variable, we have*

(a) $\partial_2 F(\bar{x}, \bar{y})(w) = P'_{\mathbb{B}}(\bar{y})(u) = u - \langle J(\bar{y}), u\rangle_X \bar{y}$, *if* $u \in \bar{y}^{\uparrow}$;
(b) $\partial_2 F(\bar{x}, \bar{y})(\bar{y}) = P'_{\mathbb{B}}(\bar{y})(\bar{y}) = \theta_X$ *(this is the special case in* (a) *with* $u = \bar{y}$);
(c) $\partial_2 F(\bar{x}, \bar{y})(w) = P'_{\mathbb{B}}(\bar{y})(u) = u$, *if* $u \in \bar{y}^{\downarrow}$.

*Proof*. This theorem follows Theorem 3.1 in [18] and Lemma 2.10. □

In order to construct another counter example by using Theorems 3.2 and 4.1 in [21], we review the concepts of the normalized duality mapping $J$: $l_p \to l_q$ and a duality mapping $\mathcal{J}$: $l_p \to l_q$, which is introduced in [21].

Throughout this paper, if not otherwise is stated, we always let $\mathbb{N}$, $\mathbb{R}$ and $\mathbb{R}_+$ respectively denote the set of nonnegative integers, the set of real numbers and the set of nonnegative real numbers. For given positive numbers $p$ and $q$ with $1 < p, q < \infty$ satisfying $\frac{1}{p} + \frac{1}{p}$ =1, let $(l_p, \|\cdot\|_p)$ be the uniformly convex and uniformly smooth Banach space with topological dual space $(l_q, \|\cdot\|_q)$, which have the same null element $\theta$ = (0, 0, …). Let $\langle\cdot,\ \cdot\rangle$ denote the pairing between $l_q$ and $l_p$. The normalized duality mapping $J$: $l_p \to l_q$ is a single-valued mapping and it satisfies that

$$\langle J(x),\ x\rangle = \|J(x)\|_q^2 = \|x\|_p^2, \text{ for any } x \in l_p.$$

Furthermore, the normalized duality mapping $J$: $l_p \to l_q$ satisfies that $J(\theta) = \theta$ and for any $x = (t_1, t_2, \ldots)$ $\in l_p\backslash\{\theta\}$, $J(x)$ is explicitly represented as follows.

$$J(x) = \left(\frac{|t_1|^{p-1}\text{sign}(t_1)}{\|x\|_p^{p-2}},\ \frac{|t_2|^{p-1}\text{sign}(t_2)}{\|x\|_p^{p-2}}, \ldots\right) = \left(\frac{|t_1|^{p-2}t_1}{\|x\|_p^{p-2}},\ \frac{|t_2|^{p-2}t_2}{\|x\|_p^{p-2}}, \ldots\right).$$

The normalized duality mapping $J$: $l_p \to l_q$ has the following differentiation properties (See [21]).

**Part (i) in Theorem 3.2 in [21]**. *Let* $1 < p, q < \infty$ *and* $\frac{1}{p} + \frac{1}{q} = 1$. *The normalized duality mapping* $J$: $l_p \to l_q$ *has the following differentiation properties.*

(i) *J is Gâteaux differentiable at* $\theta$ *and*

$$J'(\theta_p)(v) = J(v), \textit{ for any } v \in l_p\backslash\{\theta\}.$$

**Theorem 4.1 in [21]**. *Let* $1 < p, q < \infty$ *and* $\frac{1}{p} + \frac{1}{q} = 1$. *J is not Fréchet differentiable at* $\theta$; *that is,* $\nabla J(\theta)$ *does not exist.*

In [21], a duality mapping $\mathcal{J}$: $l_p \to l_q$ is defined, for any $x = (t_1, t_2, \ldots) \in l_p$, by

$$\mathcal{J}(x) = \{|t_n|^{p-1}\text{sign}(t_n)\}_{n=1}^{\infty} \in l_q.$$

The duality mapping $\mathcal{J}$: $l_p \to l_q$ has the following Fréchet differentiability (See section 6 in [21]).

**Theorem 6.1 in [21].** *Let* $3 < p < \infty$. $\mathcal{J}$ *is Fréchet differentiable on* $l_p$ *and, for any* $\bar{x} = (\bar{t}_1, \bar{t}_2, \ldots) \in l_p$,

$$\nabla\mathcal{J}(\bar{x}) = \frac{p}{q}\begin{pmatrix} |\bar{t}_1|^{p-2} & 0 & \ldots \\ 0 & |\bar{t}_2|^{p-2} & 0 \\ \vdots & 0 & \ddots \end{pmatrix}.$$

*More precisely,* $\nabla\mathcal{J}(\bar{x})$: $l_p \to l_q$ *is a pointwise multiplication operator on* $l_p$ *and, for* $x = \{t_n\}_{n=1}^{\infty} \in l_p$,

$$\nabla\mathcal{J}(\bar{x})(x) = \left\{\frac{p}{q}|\bar{t}_n|^{p-2}t_n\right\}_{n=1}^{\infty} \in l_q.$$

By Theorems 4.1 and 4.7 in [21], we construct another counter-example to show that partial Gâteaux differentiability does not imply the partial Gâteaux differentiability.

**Theorem 2.12**. *Let* $3 < p < \infty$. *Define* $F$: $l_p \times l_p \to l_q$ *as follows*:

$$F(x, y) = J(x) + \mathcal{J}(y), \text{ for any } (x, y) \in l_p \times l_p.$$

*Then, for any* $\bar{x} = (\bar{t}_1, \bar{t}_2, \dots) \in l_p$, $\bar{y} = (\bar{s}_1, \bar{s}_2, \dots) \in l_p$, *we have*

(i) (a) *F is partially Gâteaux differentiable at* $(\theta, \bar{y})$ *with respect to the first variable and*

$$\partial_1 F(\theta, \bar{y})(v) = J(v), \text{ for any } v \in l_p \backslash \{\theta\};$$

(b) *F is not partially Fréchet differentiable at* $(\theta, \bar{y})$ *with respect to the first variable and*

$$\nabla_1 F(\theta, \bar{y}) \text{ does not exist}.$$

(ii) (a) *F is partially Gâteaux differentiable on* $l_p \times l_p$ *with respect to the second variable and*

$$\partial_2 F(\bar{x}, \bar{y})(v) = \nabla_2 F(\bar{x}, \bar{y})(v) = \{(p-1)|\bar{s}_n|^{p-2}t_n\}_{n=1}^{\infty} \in l_q, \text{ for } v = \{t_n\}_{n=1}^{\infty} \in l_p \backslash \{\theta\}.$$

(b) *F is partially Fréchet differentiable on* $l_p \times l_p$ *with respect to the second variable and*

$$\nabla_2 F(\bar{x}, \bar{y}) = \nabla\mathcal{J}(\bar{y}) = \begin{pmatrix} |\bar{s}_1|^{p-2} & 0 & \cdots \\ 0 & |\bar{s}_2|^{p-2} & 0 \\ \vdots & 0 & \ddots \end{pmatrix}.$$

*Proof*. The proof of this theorem is by Theorems 3.2 and 4.1 in [21]. It is omitted here. □

## 3. Partial Gâteaux and Fréchet Derivatives in Banach Space $l_p$

### 3.1. Partial Gâteaux and Fréchet Derivatives of a Polynomial Type Operator in $l_p$

For given positive numbers $p$ and $q$ with $1 < p, q < \infty$ satisfying $\frac{1}{p} + \frac{1}{p}$ =1, the uniformly convex and uniformly smooth Banach spaces $(l_p, \|\cdot\|_p)$ and $(l_q, \|\cdot\|_q)$ are mutually topological dual space each other. Let $\langle\cdot, \cdot\rangle$ denote the real pairing between $l_q$ and $l_p$. We define a mapping $\mathcal{P}$: $l_p \to \mathbb{R}_+$ as follows

$$\mathcal{P}(x) = \sup\{|t_n|: n = 1, 2, \dots\}, \text{ for any } x = \{t_n\}_{n=1}^{\infty} \in l_p.$$

It is clear that $\mathcal{P}$ is well-defined on $l_p$ and it satisfies that,

$$0 \leq \mathcal{P}(x) < \infty, \text{ for any } x = \{t_n\}_{n=1}^{\infty} \in l_p.$$

Furthermore, by the summability of the $p^{\text{th}}$ power, for any $x = \{t_n\}_{n=1}^{\infty} \in l_p$, we have that,

$$\mathcal{P}(x) = \sup\{|t_n|: n = 1, 2, \dots\} = \max\{|t_n|: n = 1, 2, \dots\}.$$

In this section, we use Definitions 2.1 and 2.3 to find the explicit formulas of partial Gâteaux and Fréchet derivatives for some polynomial type two-variable operators in Banach space $l_p$, for $1 < p < \infty$, which are applications of partial Gâteaux and Fréchet differentiation. We start with power operators.

Let $m$ be a positive integer. We define a two-variable $m$th power operator $Q^m: l_p \times l_p \to l_p$ as follows.

$$Q^m(x, y) = \{(t_n^2 + s_n)^m\}_{n=1}^{\infty}, \text{ for any } x = \{t_n\}_{n=1}^{\infty}, y = \{s_n\}_{n=1}^{\infty} \in l_p. \tag{3.1}$$

**Theorem 3.1**. *Let m be a positive integer. Then, $Q^m$ is partially Fréchet differentiable on $l_p$ with each variable and for any $(\bar{x}, \bar{y}) \in l_p \times l_p$ with $\bar{x} = (\bar{t}_1, \bar{t}_2, \dots)$ and $\bar{y} = (\bar{s}_1, \bar{s}_2, \dots)$, we have*

(i) *The partial Fréchet derivative of $Q^m$ with respect to the first variable satisfies that,*

$$\nabla_1(Q^m)(\bar{x}, \bar{y}) = \begin{pmatrix} 2m\bar{t}_1(\bar{t}_1^2 + \bar{s}_1)^{m-1} & 0 & \dots \\ 0 & 2m\bar{t}_2(\bar{t}_2^2 + \bar{s}_2)^{m-1} & 0 \\ \vdots & 0 & \ddots \end{pmatrix}. \tag{3.2}$$

*More precisely, $\nabla_1(Q^m)(\bar{x}, \bar{y})$ is a continuous and linear operator on $l_p$ such that, for any $v = \{t_n\}_{n=1}^{\infty} \in l_p$,*

$$\nabla_1(Q^m)(\bar{x}, \bar{y})(v) = \{t_n\}_{n=1}^{\infty} \begin{pmatrix} 2m\bar{t}_1(\bar{t}_1^2 + \bar{s}_1)^{m-1} & 0 & \dots \\ 0 & 2m\bar{t}_2(\bar{t}_2^2 + \bar{s}_2)^{m-1} & 0 \\ \vdots & 0 & \ddots \end{pmatrix}$$

$$= \{2m\bar{t}_n(\bar{t}_n^2 + \bar{s}_n)^{m-1} t_n\}_{n=1}^{\infty} \in l_p. \tag{3.3}$$

(ii) *The partial Fréchet derivative of $Q^m$ with respect to the second variable satisfies that,*

$$\nabla_2(Q^m)(\bar{x}, \bar{y}) = \begin{pmatrix} m(\bar{t}_1^2 + \bar{s}_1)^{m-1} & 0 & \dots \\ 0 & m(\bar{t}_2^2 + \bar{s}_2)^{m-1} & 0 \\ \vdots & 0 & \ddots \end{pmatrix}.$$

*More precisely, $\nabla_2(Q^m)(\bar{x}, \bar{y})$ is a continuous and linear operator on $l_p$ such that, for any $u = \{s_n\}_{n=1}^{\infty} \in l_p$,*

$$\nabla_2(Q^m)(\bar{x}, \bar{y})(u) = \{s_n\}_{n=1}^{\infty} \begin{pmatrix} m(\bar{t}_1^2 + \bar{s}_1)^{m-1} & 0 & \dots \\ 0 & m(\bar{t}_2^2 + \bar{s}_2)^{m-1} & 0 \\ \vdots & 0 & \ddots \end{pmatrix}$$

$$= \{m(\bar{t}_n^2 + \bar{s}_n)^{m-1} s_n\}_{n=1}^{\infty} \in l_p. \tag{3.4}$$

*Proof*. Proof of (i). Let $\bar{x} = (\bar{t}_1, \bar{t}_2, \dots)$ and $\bar{y} = (\bar{s}_1, \bar{s}_2, \dots) \in l_p$ be arbitrarily given. Let $\nabla_1(Q^m)(\bar{x}, \bar{y})$ be represented by (3.2) with definition in (3.22). Then, for any $u = \{s_n\}_{n=1}^{\infty} \in l_p \backslash \{\theta_p\}$with $0 < \|u\|_p < 1$, we calculate

$$\frac{\|Q^m(\bar{x}+u, \bar{y}) - Q^m(\bar{x}, \bar{y}) - \nabla_1(Q^m)\,(\bar{x}, \bar{y})(u)\|_p}{\|u\|_p}$$

$$= \frac{\left\|\left\{\left((\bar{t}_n+s_n)^2+\bar{s}_n\right)^m-\left(\bar{t}_n^2+\bar{s}_n\right)^m-2m\bar{t}_n\left(\bar{t}_n^2+\bar{s}_n\right)^{m-1}s_n\right\}_{n=1}^{\infty}\right\|_p}{\|u\|_p}$$

$$= \frac{\left\|\left\{\left((\bar{t}_n^2+\bar{s}_n)+(2\bar{t}_n s_n+s_n^2)\right)^m-\left(\bar{t}_n^2+\bar{s}_n\right)^m-2m\bar{t}_n\left(\bar{t}_n^2+\bar{s}_n\right)^{m-1}s_n\right\}_{n=1}^{\infty}\right\|_p}{\|u\|_p}$$

$$= \frac{\left\|\left\{\left(\bar{t}_n^2+\bar{s}_n\right)^m+\binom{m}{1}\left(\bar{t}_n^2+\bar{s}_n\right)^{m-1}(2\bar{t}_n s_n+s_n^2)+\binom{m}{2}\left(\bar{t}_n^2+\bar{s}_n\right)^{m-2}\left(2\bar{t}_n s_n+s_n^2\right)^2+\cdots+\binom{m}{m}(2\bar{t}_n s_n+s_n^2)^m\right.}{\|u\|_p}$$

$$\frac{\left.-\left(\bar{t}_n^2+\bar{s}_n\right)^m-2m\bar{t}_n\left(\bar{t}_n^2+\bar{s}_n\right)^{m-1}s_n\right\}_{n=1}^{\infty}\Big\|_p}{\|u\|_p}$$

$$= \frac{\left\|\left\{\binom{m}{1}\left(\bar{t}_n^2+\bar{s}_n\right)^{m-1}(2\bar{t}_n s_n+s_n^2)+\binom{m}{2}\left(\bar{t}_n^2+\bar{s}_n\right)^{m-2}(2\bar{t}_n s_n+s_n^2)^2+\cdots+\binom{m}{m}(2\bar{t}_n s_n+s_n^2)^m-2m\bar{t}_n\left(\bar{t}_n^2+\bar{s}_n\right)^{m-1}s_n\right\}_{n=1}^{\infty}\right\|_p}{\|u\|_p}$$

$$= \frac{\left\|\left\{m\left(\bar{t}_n^2+\bar{s}_n\right)^{m-1}s_n^2+\binom{m}{2}\left(\bar{t}_n^2+\bar{s}_n\right)^{m-2}(2\bar{t}_n s_n+s_n^2)^2+\cdots+\binom{m}{m}(2\bar{t}_n s_n+s_n^2)^m\right\}_{n=1}^{\infty}\right\|_p}{\|u\|_p}$$

$$= \frac{\left\|\left\{s_n^2\left(m\left(\bar{t}_n^2+\bar{s}_n\right)^{m-1}+\binom{m}{2}\left(\bar{t}_n^2+\bar{s}_n\right)^{m-2}(2\bar{t}_n+s_n)^2+\cdots+\binom{m}{m}(2\bar{t}_n+s_n)^m s_n^{m-2}\right)\right\}_{n=1}^{\infty}\right\|_p}{\|u\|_p}$$

$$= \frac{\left\|\left\{s_n^2\left(m\left(\bar{t}_n^2+\bar{s}_n\right)^{m-1}+\binom{m}{2}\left(\bar{t}_n^2+\bar{s}_n\right)^{m-2}(2\bar{t}_n+s_n)^2+\cdots+\binom{m}{m}(2\bar{t}_n+s_n)^m s_n^{m-2}\right)\right\}_{n=1}^{\infty}\right\|_p}{\|u\|_p}$$

$$= \frac{\left(\sum_{n=1}^{\infty}\left|s_n^2\left(m\left(\bar{t}_n^2+\bar{s}_n\right)^{m-1}+\binom{m}{2}\left(\bar{t}_n^2+\bar{s}_n\right)^{m-2}(2\bar{t}_n+s_n)^2+\cdots+\binom{m}{m}(2\bar{t}_n+s_n)^m s_n^{m-2}\right)\right|^p\right)^{\frac{1}{p}}}{\|u\|_p}$$

$$\leq \frac{\left(\sum_{n=1}^{\infty}|s_n|^{2p}\left(m\left|\bar{t}_n^2+\bar{s}_n\right|^{m-1}+\binom{m}{2}\left|\bar{t}_n^2+\bar{s}_n\right|^{m-2}(2|\bar{t}_n|+1)^2+\cdots+\binom{m}{m}(2|\bar{t}_n|+1)^m\right)^p\right)^{\frac{1}{p}}}{\|u\|_p}$$

$$\leq \frac{\left(\sum_{n=1}^{\infty}|s_n|^{2p}\left(m(\mathcal{P}(\bar{x}^2)+\mathcal{P}(\bar{y}))^{m-1}+\binom{m}{2}(\mathcal{P}(\bar{x}^2)+\mathcal{P}(\bar{y}))^{m-2}(2\mathcal{P}(\bar{x})+1)^2+\cdots+\binom{m}{m}(2\mathcal{P}(\bar{x})+1)^m\right)^p\right)^{\frac{1}{p}}}{\|u\|_p}$$

$$\leq \frac{\left(\sum_{n=1}^{\infty}|s_n|^{2p}\left((\mathcal{P}(\bar{x}^2)+2\mathcal{P}(\bar{x})+\mathcal{P}(\bar{y})+1)^m\right)^p\right)^{\frac{1}{p}}}{\|u\|_p}$$

$$= \frac{(\mathcal{P}(\bar{x}^2)+2\mathcal{P}(\bar{x})+\mathcal{P}(\bar{y})+1)^m\left(\sum_{n=1}^{\infty}|s_n|^{2p}\right)^{\frac{1}{p}}}{\|u\|_p}$$

$$= \frac{\left(\mathcal{P}(\bar{x}^2)+2\mathcal{P}(\bar{x})+\mathcal{P}(\bar{y})+1\right)^m\left(\sum_{n=1}^{\infty}|s_n|^p|s_n|^p\right)^{\frac{1}{p}}}{\|u\|_p}$$

$$\leq \frac{\left(\mathcal{P}(\bar{x}^2)+2\mathcal{P}(\bar{x})+\mathcal{P}(\bar{y})+1\right)^m \left(\sum_{n=1}^{\infty}|s_n|^p\left(\sum_{i=1}^{\infty}|s_i|^p\right)\right)^{\frac{1}{p}}}{\|u\|_p}$$

$$= \frac{\left(\mathcal{P}(\bar{x}^2)+2\mathcal{P}(\bar{x})+\mathcal{P}(\bar{y})+1\right)^m \left(\sum_{n=1}^{\infty}|s_n|^p\right)^{\frac{1}{p}}\left(\left(\sum_{i=1}^{\infty}|s_i|^p\right)\right)^{\frac{1}{p}}}{\|u\|_p}$$

$$= (\mathcal{P}(\bar{x}^2) + 2\mathcal{P}(\bar{x}) + \mathcal{P}(\bar{y}) + 1)^m\|u\|_p.$$

This implies that

$$\lim_{u\to\theta_X} \frac{\|\mathcal{Q}^m(\bar{x}+u,\bar{y})-\mathcal{Q}^m(\bar{x},\bar{y})-\nabla_1(\mathcal{Q}^m)\ (\bar{x},\bar{y})(u)\|_p}{\|u\|_p} = 0.$$

Part (ii) can be similarly proved. □

By the connection between partial Gâteaux and Fréchet derivatives given in (2.10) in Theorem 2.5, we have the following corollary of Theorem 3.1.

**Corollary 3.2**. *Let m be a positive integer. Then, $\mathcal{Q}^m$ is partially Gâteaux differentiable on $l_p$ with each variable and, for any $(\bar{x},\bar{y}) \in l_p \times l_p$ with $\bar{x} = (\bar{t}_1, \bar{t}_2, ...)$ and $\bar{y} = (\bar{s}_1, \bar{s}_2, ...)$, we have*

(i) *The partial Gâteaux derivative of $\mathcal{Q}^m$ with respect to the first variable satisfies that, for any $v = \{t_n\}_{n=1}^{\infty} \in l_p\backslash\{\theta_p\}$, we have*

$$\partial_1\mathcal{Q}^m(\bar{x},\bar{y})(v) = \nabla_1(\mathcal{Q}^m)(\bar{x},\bar{y})(v) = \{t_n\}_{n=1}^{\infty}\begin{pmatrix} 2m\bar{t}_1(\bar{t}_1^2+\bar{s}_1)^{m-1} & 0 & \cdots \\ 0 & 2m\bar{t}_2(\bar{t}_2^2+\bar{s}_2)^{m-1} & 0 \\ \vdots & 0 & \ddots \end{pmatrix}.$$

(ii) *The partial Gâteaux derivative of $\mathcal{Q}^m$ with respect to the second variable satisfies that, for $u = \{s_n\}_{n=1}^{\infty} \in l_p\backslash\{\theta_p\}$,*

$$\partial_2\mathcal{Q}^m(\bar{x},\bar{y})(u) = \nabla_2(\mathcal{Q}^m)(\bar{x},\bar{y})(u) = \{s_n\}_{n=1}^{\infty}\begin{pmatrix} m(\bar{t}_1^2+\bar{s}_1)^{m-1} & 0 & \cdots \\ 0 & m(\bar{t}_2^2+\bar{s}_2)^{m-1} & 0 \\ \vdots & 0 & \ddots \end{pmatrix}.$$

Let *m* be a positive integer. Let $a_m, a_{m-1}, \dots a_1$ be real numbers. By the power operators defined in (3.1), define a polynomial type operator $\mathcal{Q}_m: l_p \times l_p \to l_p$, for $x = \{t_n\}_{n=1}^{\infty}$, $y = \{s_n\}_{n=1}^{\infty} \in l_p$, by

$$\mathcal{Q}_m(x,y) = a_m\mathcal{Q}^m(x,y) + a_{m-1}\mathcal{Q}^{m-1}(x,y) + \ldots + a_1\mathcal{Q}(x,y) = \left\{\sum_{i=1}^{m} a_i(t_n^2+s_n)^i\right\}_{n=1}^{\infty}. \quad (3.5)$$

By the linearity property of partial Fréchet derivative given in (2.9) in Theorem 2.4 and by (3.2) in Theorem 3.1, we obtain the following results.

**Corollary 3.3**. *Let m be a positive integer and let $\mathcal{Q}_m$ be defined by* (3.5). *Then, $\mathcal{Q}_m$ is partially Fréchet differentiable on $l_p$ with each variable such that for given $(\bar{x},\bar{y}) \in l_p \times l_p$ with $\bar{x} = (\bar{t}_1, \bar{t}_2, ...)$ and $\bar{y} = (\bar{s}_1, \bar{s}_2, ...)$, we have*

(i) *The partial Fréchet derivative of $\mathcal{Q}_m$ with respect to the first variable $\nabla_1(\mathcal{Q}_m)(\bar{x},\bar{y}): l_p \to l_p$ satisfies that,*

$$\nabla_1(Q_m)(\bar{x},\bar{y}) = 2\begin{pmatrix} \sum_{i=1}^{m} a_i i\, \bar{t}_1(\bar{t}_1^2+\bar{s}_1)^{i-1} & 0 & \dots \\ 0 & \sum_{i=1}^{m} a_i i\, \bar{t}_2(\bar{t}_2^2+\bar{s}_2)^{i-1} & 0 \\ \vdots & 0 & \ddots \end{pmatrix}. \quad (3.6)$$

(ii) *The partial Fréchet derivative of $Q_m$ with respect to the second variable $\nabla_2(Q_m)(\bar{x},\bar{y}): l_p \to l_p$ satisfies that,*

$$\nabla_2(Q_m)(\bar{x},\bar{y}) = \begin{pmatrix} \sum_{i=1}^{m} a_i i\, (\bar{t}_1^2+\bar{s}_1)^{i-1} & 0 & \dots \\ 0 & \sum_{i=1}^{m} i a_i\, (\bar{t}_2^2+\bar{s}_2)^{i-1} & 0 \\ \vdots & 0 & \ddots \end{pmatrix}. \quad (3.7)$$

Similar to Corollary 3.2, by the connection between partial Gâteaux and Fréchet derivatives given in (2.10) in Theorem 2.5, we have the following corollary of Corollary 3.3.

**Corollary 3.4**. *Let m be a positive integer. Then, $Q_m$ is partially Gâteaux differentiable on $l_p$ with each variable such that for given $(\bar{x},\bar{y}) \in l_p \times l_p$ with $\bar{x} = (\bar{t}_1, \bar{t}_2, \dots)$ and $\bar{y} = (\bar{s}_1, \bar{s}_2, \dots)$, the partial Gâteaux derivatives of $Q_m$ at $(\bar{x},\bar{y})$ have the following properties*:

(i) *The partial Gâteaux derivative of $Q_m$ with respect to the first variable satisfies that, for any $v = \{t_n\}_{n=1}^{\infty} \in l_p$,*

$$\partial_1 Q_m(\bar{x},\bar{y})(v) = \nabla_1(Q_m)(\bar{x},\bar{y})(v) = 2\left\{\sum_{i=1}^{m} a_i i\, \bar{t}_n(\bar{t}_n^2+\bar{s}_n)^{i-1} t_n\right\}_{n=1}^{\infty} \in l_p. \quad (3.8)$$

(ii) *The partial Gâteaux derivative of $Q_m$ with respect to the second variable satisfies that, for any $u = \{s_n\}_{n=1}^{\infty} \in l_p$,*

$$\partial_2 Q_m(\bar{x},\bar{y})(u) = \nabla_2(Q_m)(\bar{x},\bar{y})(u) = \left\{\sum_{i=1}^{m} i a_i\, (\bar{t}_n^2+\bar{s}_n)^{i-1} s_n\right\}_{n=1}^{\infty} \in l_p. \quad (3.9)$$

## 4. Partially Ordered extrema and Partially Ordered Optimization Problems

In this section, we will introduce the concepts of partially ordered extrema of multivariable and single-valued mappings between partially ordered Banach spaces. We consider some applications of partial differentiation to partially ordered optimization problems, which are equivalent to ordered variational inequality problems with multivariable mappings in partially ordered Banach spaces.

By following the notations in section 2, let $(X, \|\cdot\|_X)$ and $(Y, \|\cdot\|_Y)$ be Banach spaces and let $(Z, \|\cdot\|_Z, \preccurlyeq_M)$ be a partially ordered Banach space, in which the partial order $\preccurlyeq_M$ is generated by a nonempty closed convex and pointed cone $M$ in $Z$. For given $x_0 \in X$ and $r > 0$, let $B_X(x_0, r)$ denote the closed ball in $X$ with radius $r$ and centered at $x_0$. For given $y_0 \in Y$, $z_0 \in Z$ and $r > 0$, the closed balls $B_Y(y_0, r)$ and $B_Z(z_0, r)$ with radius $r$ in $Y$ and $Z$ are similarly defined, respectively. Let $A$ and $B$ be nonempty convex open subsets in $X$ and $Y$ respectively.

**Definition 4.1**. Let $F: A \times B \to Z$ be a single-valued mapping. Let $(\bar{x},\bar{y}) \in A \times B$. If

$$F(\bar{x},\bar{y}) \preccurlyeq_M F(x,\bar{y}), \text{ for every } x \in A, \quad (4.1)$$

then, $F$ is said to have a partial ordered-minimum value at $(\bar{x},\bar{y})$ with respect to the first variable on $A \times B$. A partial ordered-maximum value of $F$ with respect to the first variable on $A \times B$ can be similarly defined. The set of all points at which $F$ has a partial ordered-minimum value (partial ordered-maximum value) at $(\bar{x}, \bar{y})$ with respect to the first variable on $A \times B$ is denoted by $E_1^{\wedge}(F, A, B)$ ($E_1^{\vee}(F, A, B)$). In particular, $E_1^{\wedge}(F, X, Y)$ and $E_1^{\vee}(F, X, Y)$ are simply denoted by $E_1^{\wedge}(F)$ and $E_1^{\vee}(F)$, respectively.

If $F$ has either partial ordered-minimum value or partial ordered-maximum value at $(\bar{x}, \bar{y})$ with respect to the first variable on $A \times B$, then $F$ is said to have a partial ordered-extreme value at $(\bar{x}, \bar{y})$ with respect to the first variable on $A \times B$, and such a point $(\bar{x}, \bar{y})$ is called a partial ordered-extreme point of $F$ with respect to the first variable on $A \times B$. The set of all partial ordered-extreme points of $F$ with respect to the first variable on $A \times B$ is denoted by $E_1(F, A, B)$. In particular, $E_1(F, X, Y)$ is simply denoted by $E_1(F)$. The following properties are satisfied.

$$E_1(F, A, B) = E_1^{\wedge}(F, A, B) \cup E_1^{\vee}(F, A, B) \quad \text{and} \quad E_1(F) = E_1^{\wedge}(F) \cup E_1^{\vee}(F).$$

Similarly to (4.1), if the point $(\bar{x}, \bar{y}) \in A \times B$ satisfies

$$F(\bar{x}, \bar{y}) \preccurlyeq_M F(\bar{x}, y), \text{ for every } y \in B. \tag{4.2}$$

then, $F$ is said to have a partial ordered-minimum value at $(\bar{x}, \bar{y})$ with respect to the second variable on $A \times B$. A partial ordered-maximum value of $F$ with respect to the second variable on $A \times B$ can be similarly defined. Then, with respect to the second variable, $E_2^{\wedge}(F, A, B)$, $E_2^{\vee}(F, A, B)$, $E_2^{\wedge}(F)$ and $E_2^{\vee}(F)$ are similarly defined.

If $F$ has either partial ordered-minimum value or partial ordered-maximum value at $(\bar{x}, \bar{y})$ with respect to the second variable on $A \times B$, then $F$ is said to have a partial ordered-extreme value at $(\bar{x}, \bar{y})$ with respect to the second variable on $A \times B$, and such a point $(\bar{x}, \bar{y})$ is called a partial ordered-extreme point of $F$ with respect to the second variable on $A \times B$. The set of all partial ordered-extreme points of $F$ with respect to the second variable on $A \times B$ is denoted by $E_2(F, A, B)$. In particular, $E_2(F, X, Y)$ is simply denoted by $E_2(F)$. The following properties are satisfied.

$$E_2(F, A, B) = E_2^{\wedge}(F, A, B) \cup E_2^{\vee}(F, A, B) \quad \text{and} \quad E_2(F) = E_2^{\wedge}(F) \cup E_2^{\vee}(F)$$

More strictly, for the given mapping $F$, if $(\bar{x}, \bar{y}) \in A \times B$ satisfies the following condition

$$F(\bar{x}, \bar{y}) \preccurlyeq_M F(x, y), \text{ for every } x \in A \text{ and for every } y \in B,$$

then $F$ is said to have an ordered-minimum value at point $(\bar{x}, \bar{y})$. The ordered-maximum value of $F$ can be similarly defined. At such a point, for either ordered-minimum, or ordered-maximum, $F$ is said to have an ordered-extremum. The set of all ordered-extreme points of $F$ on $A \times B$ is denoted by $E(F, A, B)$. In particular, $E(F, X, Y)$ is simply denoted by $E(F)$. One can see that, with respect to given $(\bar{x}, \bar{y}) \in A \times B$

$$E(F, A, B) \subseteq E_1(F, A, B) \cap E_2(F, A, B) \quad \text{and} \quad E(F) \subseteq E_1(F) \cap E_2(F). \tag{4.3}$$

Furthermore, suppose that $F$ is partially Gâteaux differentiable with respect to the first (second) variable, in Theorem 4.3, we will prove the properties of partially ordered extrema, which actually extend the properties of critical points in calculus.

$F$ has a partial ordered-extreme value at point $(\bar{x}, \bar{y})$ with respect to the first (second) variable

$\Longrightarrow$ $(\bar{x}, \bar{y})$ is a generalized partially critical point of $F$ with respect to the first (second) variable. (4.4)

**Theorem 4.3**. *Let $F$: $A \times B \to Z$ be a single-valued mapping. Let $(\bar{x}, \bar{y}) \in A \times B$.*

(i) *If $F$ is partially Gâteaux differentiable with respect to the first variable at point $(\bar{x}, \bar{y})$, then*

$$(\bar{x}, \bar{y}) \in E_1(F, A, B) \quad \Longrightarrow \quad (\bar{x}, \bar{y}) \in \partial_1^C F.$$

(ii) *If F is partially Gâteaux differentiable with respect to the second variable at point* $(\bar{x}, \bar{y})$, *then*

$$(\bar{x}, \bar{y}) \in E_2(F, A, B) \quad \Longrightarrow \quad (\bar{x}, \bar{y}) \in \partial_2^C F.$$

*Proof*. Proof of (i). Suppose that $F$ is partially Gâteaux differentiable with respect to the first variable at point $(\bar{x}, \bar{y})$ and $F$ has a partial ordered-extreme at $(\bar{x}, \bar{y})$ with respect to the first variable. By definition, it follows that $F$ has either a partial ordered-minimum value or a partial ordered-maximum value at point $(\bar{x}, \bar{y})$ with respect to the first variable. Suppose that $F$ has a partial ordered-minimum value at point $(\bar{x}, \bar{y})$ with respect to the first variable. By which, we want to prove

$$\partial_1 F(\bar{x}, \bar{y})(v) = \theta_Z, \text{ for any } v \in X \backslash \{\theta_X\}.$$

At first, we prove

$$\partial_1 F(\bar{x}, \bar{y})(v) \succcurlyeq_M \theta_Z, \text{ for any } v \in X \backslash \{\theta_X\}. \tag{4.5}$$

Assume, by the way of contradiction, that there is $v \in X \backslash \{\theta_X\}$ such that $\frac{\partial F(x,y)}{\partial x}(\bar{x}, \bar{y})(v) \not\succcurlyeq_M \theta_Z$. This is that $\partial_1 F(\bar{x}, \bar{y})(v) \notin M$, this implies $\partial_1 F(\bar{x}, \bar{y})(v) \in Z \backslash M$. Since $Z \backslash M$ is an open subset of $Z$, by $\partial_1 F(\bar{x}, \bar{y})(v) \in Z \backslash M$, there is $\varepsilon > 0$ such that, $B_Z^o(\partial_1 F(\bar{x}, \bar{y})(v), \varepsilon) \subseteq Z \backslash M$, which defines the topological interior of the closed ball $B_Z(\partial_1 F(\bar{x}, \bar{y})(v), \varepsilon)$ in $Z$ with radius $\varepsilon$ and centered at point $\partial_1 F(\bar{x}, \bar{y})(v)$. More precisely, we have

$$B_Z^o(\partial_1 F(\bar{x}, \bar{y})(v), \varepsilon) = \{z \in Z : \|z - \partial_1 F(\bar{x}, \bar{y})(v)\|_Z < \varepsilon\}. \tag{4.6}$$

The condition $B_Z^o(\partial_1 F(\bar{x}, \bar{y})(v), \varepsilon) \subseteq Z \backslash M$ and (4.6) induce

$$B_Z^o(\partial_1 F(\bar{x}, \bar{y})(v), \varepsilon) \cap M = \emptyset. \tag{4.7}$$

(4.6) and (4.7) together imply that, for $z \in Z$, we have

$$z \in M \quad \Longrightarrow \quad \|z - \partial_1 F(\bar{x}, \bar{y})(v)\|_Z \geq \varepsilon. \tag{4.8}$$

The condition that $F$ is partially Gâteaux differentiable with respect to the first variable at point $(\bar{x}, \bar{y})$ induces that $F$ is partially Gâteaux directionally differentiable with respect to the first variable at point $(\bar{x}, \bar{y})$ along the given (fixed) direction $v$. By definition and by the condition that $A$ is an open subset in $X$ and $\bar{x} \in A$, for the above given $\varepsilon > 0$ and with the fixed $\bar{y} \in Y$, there is $\delta > 0$ such that, for real $t$, we have

$$0 < |t| < \delta \quad \Longrightarrow \bar{x} + tv \in A \text{ and } \left\| \frac{F(\bar{x}+tv,\bar{y}) - F(\bar{x},\bar{y})}{t} - \partial_1 F(\bar{x}, \bar{y})(v) \right\|_Z < \varepsilon. \tag{4.9}$$

In particular, let $t > 0$ in (4.9), we obtain that

$$0 < t < \delta \quad \Longrightarrow \quad \bar{x} + tv \in A \text{ and } \left\| \frac{F(\bar{x}+tv,\bar{y}) - F(\bar{x},\bar{y})}{t} - \partial_1 F(\bar{x}, \bar{y})(v) \right\|_Z < \varepsilon. \tag{4.10}$$

By the assumption that $F$ has a partial ordered-minimum value at point $(\bar{x}, \bar{y})$ with respect to the first variable, we must have $F(\bar{x} + tv, \bar{y}) - F(\bar{x}, \bar{y}) \succcurlyeq_M \theta_Z$ and therefore,

$$0 < t < \delta \quad \Longrightarrow \quad \frac{F(\bar{x}+tv,\bar{y}) - F(\bar{x},\bar{y})}{t} \in M. \tag{4.11}$$

For $t > 0$, by (4.11), this creates a contradiction between (4.10) and (4.8). Hence, we must have

$$\partial_1 F(\bar{x},\bar{y})(v) \in M \quad \text{and} \quad \partial_1 F(\bar{x},\bar{y})(v) \succcurlyeq_M \theta_Z.$$

This proves (4.5). Next, we prove

$$\partial_1 F(\bar{x},\bar{y})(v) \preccurlyeq_M \theta_Z, \text{ for any } v \in X\backslash\{\theta_X\}. \tag{4.12}$$

(4.12) is equivalent to

$$\partial_1 F(\bar{x},\bar{y})(v) \in -M, \text{ for any } v \in X\backslash\{\theta_X\}. \tag{4.13}$$

Assume, by the way of contradiction, that there is $v \in X\backslash\{\theta_X\}$ such that $\partial_1 F(\bar{x},\bar{y})(v) \not\preccurlyeq_M \theta_Z$. This is that $\partial_1 F(\bar{x},\bar{y})(v) \notin -M$, this implies $\partial_1 F(\bar{x},\bar{y})(v) \in Z\backslash(-M)$. Notice that $-M$ is also a closed convex and pointed cone in $Z$. Since $Z\backslash(-M)$ is an open subset of $Z$, by $\frac{\partial F(x,y)}{\partial x}(\bar{x},\bar{y})(v) \in Z\backslash(-M)$, there is $\varepsilon > 0$ such that, $B_Z^o(\partial_1 F(\bar{x},\bar{y})(v),\varepsilon) \subseteq Z\backslash(-M)$. More precisely, we have

$$B_Z^o(\partial_1 F(\bar{x},\bar{y})(v),\varepsilon) = \{z \in Z : \|z - \partial_1 F(\bar{x},\bar{y})(v)\|_Z < \varepsilon\}. \tag{4.14}$$

The condition $B_Z^o(\partial_1 F(\bar{x},\bar{y})(v),\varepsilon) \subseteq Z\backslash(-M)$ and (4.14) induce

$$B_Z^o(\partial_1 F(\bar{x},\bar{y})(v),\varepsilon) \cap (-M) = \emptyset. \tag{4.15}$$

(4.14) and (4.15) together imply that, for $z \in Z$, we have

$$z \in -M \quad \Longrightarrow \quad \|z - \partial_1 F(\bar{x},\bar{y})(v)\|_Z \geq \varepsilon. \tag{4.16}$$

The condition that $F$ is partially Gâteaux differentiable with respect to the first variable at point $(\bar{x},\bar{y})$ induces that $F$ is partially Gâteaux directionally differentiable with respect to the first variable at point $(\bar{x},\bar{y})$ along the given (fixed) direction $v$. By definition and by the condition that $A$ is an open subset in $X$ and $\bar{x} \in A$, for the above given $\varepsilon > 0$ and with the fixed $\bar{y} \in Y$, there is $\delta > 0$ such that, for real $t$, we have

$$0 < |t| < \delta \quad \Longrightarrow \bar{x} + tv \in A \ \text{ and } \ \left\|\frac{F(\bar{x}+tv,\bar{y})-F(\bar{x},\bar{y})}{t} - \partial_1 F(\bar{x},\bar{y})(v)\right\|_Z < \varepsilon. \tag{4.17}$$

In particular, let $t < 0$ in (4.17), we obtain that

$$0 < -t < \delta \quad \Longrightarrow \quad \bar{x} + tv \in A \ \text{ and } \ \left\|\frac{F(\bar{x}+tv,\bar{y})-F(\bar{x},\bar{y})}{t} - \partial_1 F(\bar{x},\bar{y})(v)\right\|_Z < \varepsilon. \tag{4.18}$$

By the assumption that $F$ has a partial ordered-minimum value at point $(\bar{x},\bar{y})$ with respect to the first variable, we must have $F(\bar{x}+tv,\bar{y}) - F(\bar{x},\bar{y}) \succcurlyeq_M \theta_Z$ and therefore,

$$0 < -t < \delta \quad \Longrightarrow \quad \frac{F(\bar{x}+tv,\bar{y})-F(\bar{x},\bar{y})}{t} \in -M. \tag{4.19}$$

For $t < 0$, by (4.19), this creates a contradiction between (4.18) and (4.16). Hence, we must have

$$\partial_1 F(\bar{x},\bar{y})(v) \in -M \quad \text{and} \quad \partial_1 F(\bar{x},\bar{y})(v) \preccurlyeq_M \theta_Z.$$

This proves (4.12). By (4.5) and (4.12), we proved that

$F$ has a partial ordered-minimum value at $(\bar{x},\bar{y})$ with respect to the first variable

$$\Longrightarrow \quad \partial_1 F(\bar{x}, \bar{y})(v) = \theta_Z, \text{ for any v} \in X\backslash\{\theta_X\}.$$

Similarly to the above property, one can prove that

$F$ has a partial ordered-maximum value at $(\bar{x}, \bar{y})$ with respect to the first variable

$$\Longrightarrow \quad \partial_1 F(\bar{x}, \bar{y})(v) = \theta_Z, \text{ for any v} \in X\backslash\{\theta_X\}.$$

Then, part (i) is proved immediately. The proof of (ii) is similar to the proof of (i) and omitted here. □

By the connection between partial Gâteaux and Fréchet derivatives given in (2.10) in Theorem 2.5, and by (4.3), we have the following corollary of Theorem 4.3.

**Corollary 4.4**. *Let $F: A\times B \to Z$ be a single-valued mapping. Let $(\bar{x}, \bar{y}) \in A \times B$.*

(i) *If F is partial Fréchet differentiable with respect to the first variable at point $(\bar{x}, \bar{y})$, then*

$$(\bar{x}, \bar{y}) \in E_1(F, A, B) \;\; \Longrightarrow \;\; (\bar{x}, \bar{y}) \in \nabla_1^C F.$$

(ii) *If F is partial Fréchet differentiable with respect to the second variable at point $(\bar{x}, \bar{y})$, then*

$$(\bar{x}, \bar{y}) \in E_2(F, A, B) \quad \Longrightarrow \quad (\bar{x}, \bar{y}) \in \nabla_2^C F$$

(iii) *If F is partial Fréchet differentiable with respect to each variable at point $(\bar{x}, \bar{y})$, then*

$$(\bar{x}, \bar{y}) \in E(F; A, B) \quad \Longrightarrow \quad (\bar{x}, \bar{y}) \in \nabla_1^C F \cap \nabla_2^C F.$$

As a corollary of Theorem 4.3 and Corollary 4.4, next proposition explicitly shows the connection between the solutions of partially ordered variational inequality problems and generalized partially critical points of single-valued mappings.

**Proposition 4.5**. *Let $F: A\times B \to Z$ be a single-valued mapping.*

(i) *Suppose that F is partial Gâteaux differentiable on $A\times B$ with respect to every variable. Then*

$$E_1(F, A, B) \subseteq \partial_1^C F \quad \text{and} \quad E_2(F, A, B) \subseteq \partial_2^C F. \tag{4.20}$$

(ii) *Suppose that F is partial Fréchet differentiable on $A\times B$ with respect to every variable. Then*

$$E_1(F, A, B) \subseteq \nabla_1^C F \quad \text{and} \quad E_2(F, A, B) \subseteq \nabla_2^C F. \tag{4.21}$$

(iii) *If F is partial Gâteaux differentiable with respect to each variable at point $(\bar{x}, \bar{y})$, then*

$$E(F; A, B) \subseteq \partial_1^C F \cap \partial_2^C F.$$

(iv) *If F is partial Fréchet differentiable on $A\times B$ with respect to every variable. Then*

$$E(F; A, B) \subseteq \nabla_1^C F \cap \nabla_2^C F.$$

*Proof*. By Theorem 4.3, Corollary 4.4 and by Notation 2.8, this proposition is immediately proved. □

## 5. Ordered Optimization of the Power Operator $Q^m$ in $l_p$

In this section, by using the Banach space $l_p$, we provider a counter example to show that the inverse of Theorem 4.3 and Proposition 4.5 do not hold. More precisely speaking, we will show that the inverse of (4.4) is not true.

In particular, let $(X, \|\cdot\|_X, \preccurlyeq_C) = (Y, \|\cdot\|_Y, \preccurlyeq_K) = (Z, \|\cdot\|_Z, \preccurlyeq_M) = (l_p, \|\cdot\|_p, \preccurlyeq_K)$ with $1 < p < \infty$ studied in Theorem 3.1, in which $K$ is the positive cone in $l_p$, by which we construct the following example.

**Example 5.1**. In this example, we consider a two-variable operator $Q^3$: $l_p \times l_p \to l_p$ as defined in (3.1). In order to show that the inverse of (4.4) is not true, we find a point $(\bar{x}, \bar{y}) \in l_p \times l_p$ that satisfies

$$(\bar{x}, \bar{y}) \in \nabla_1^C(Q^3) \quad \text{and} \quad (\bar{x}, \bar{y}) \notin E_1(Q^3).$$

Recall that let $m = 3$ in Definition (3.1), we have

$$Q^3(x, y) = \{(t_n^2 + s_n)^3\}_{n=1}^{\infty}, \text{ for any } x = \{t_n\}_{n=1}^{\infty}, y = \{s_n\}_{n=1}^{\infty} \in l_p.$$

To this end, let $\bar{x} = \left\{\frac{1}{n}\right\}_{n=1}^{\infty}$ and $\bar{y} = \left\{-\frac{1}{n^2}\right\}_{n=1}^{\infty} \in l_p$. By part (i) in Theorem 3.1, we have

$$\partial_1(Q^3)(\bar{x}, \bar{y})(x) = \nabla_1(Q^3)\,(\bar{x}, \bar{y})(x) = \left\{\frac{6}{n}0^{3-1}t_n\right\}_{n=1}^{\infty} = \theta \in l_p, \text{ for any } x = \{t_n\}_{n=1}^{\infty} \in l_p.$$

We obtain that

$$\nabla_1(Q^3)\,(\bar{x}, \bar{y}) = \theta_{XZ}.$$

This implies that the point $(\bar{x}, \bar{y})$ is a partially generalized critical point of $Q^3$ with respect to the first variable. That is,

$$(\bar{x}, \bar{y}) \in \nabla_1^C(Q^3). \tag{5.1}$$

For real number $t$, we have

$$Q^3(\bar{x}(1+t), \bar{y}) = \left\{\left(\left(\frac{1+t}{n}\right)^2 - \frac{1}{n^2}\right)^3\right\}_{n=1}^{\infty} = \left\{\left(\frac{2t+t^2}{n}\right)^3\right\}_{n=1}^{\infty}.$$

This implies that

$$t > 0 \implies Q^3(\bar{x}(1+t), \bar{y}) \succ_K \theta,$$

and

$$-2 < t < 0 \implies Q^3(\bar{x}(1+t), \bar{y}) \prec_K \theta.$$

This implies that $F$ does not have a partial ordered-extreme value at $(\bar{x}, \bar{y})$ with respect to the first variable. That is,

$$(\bar{x}, \bar{y}) \notin E_1(Q^3). \tag{5.2}$$

Hence, by (5.1) and (5.2), we obtain that the point $(\bar{x}, \bar{y})$ satisfies

$$(\bar{x}, \bar{y}) \in \nabla_1^C(Q^3) \quad \text{and} \quad (\bar{x}, \bar{y}) \in E_1(Q^3).$$

This counter example shows that the inverse of (4.4) does not hold.

Next, for the power operator $Q^m: l_p \times l_p \to l_p$ defined in (3.1), in the following theorem, we will find the explicit solutions for $E_1(Q^m)$ and $E_2(Q^m)$, which verify (4.20) and (4.21) in Proposition 4.7. Recall that by Definition 4.1, $E_1(Q^m, l_p, l_p)$ and $E_2(Q^m, l_p, l_p)$ are written as $E_1(Q^m)$ and $E_2(Q^m)$, respectively.

**Theorem 5.2**. *Let* $1 < p < \infty$. *Let* $(l_p, \|\cdot\|_p, \preccurlyeq_K)$ *be the partially ordered Banach space, in which K is the positive cone in* $l_p$. *Let m be a positive integer and let the two-variable operator* $Q^m: l_p \times l_p \to l_p$ *be defined in* (3.1).

$$Q^m(x,y) = \{(t_n^2 + s_n)^m\}_{n=1}^{\infty}, \text{ for any } x = \{t_n\}_{n=1}^{\infty}, y = \{s_n\}_{n=1}^{\infty} \in l_p. \quad (5.3)$$

(i) $m = 1$. *We have*

$$\partial_1^C Q = \nabla_1^C(Q) = \{\theta\} \times l_p, \quad \partial_2^C Q = \nabla_2^C(Q) = \emptyset; \quad (5.4)$$

*and*

$$E_1(Q) = \{\theta\} \times l_p, \quad E_2(Q) = \emptyset. \quad (5.5)$$

(ii) $m > 1$ *and m is odd. Then,*

(a) *The generalized critical points of* $Q^m$ *satisfy*

$$\partial_1^C(Q^m) = \nabla_1^C(Q^m) = (\{\theta\} \times l_p) \cup \{(\bar{x},\bar{y}) \in l_p \times l_p : \bar{y} = -\bar{x}^2\}, \quad (5.6)$$

*and*

$$\partial_2^C(Q^m) = \nabla_2^C(Q^m) = \{(\bar{x},\bar{y}) \in l_p \times l_p : \bar{y} = -\bar{x}^2\}. \quad (5.7)$$

(b) *The solution sets to the ordered variational inequality problems satisfy*

$$E_1(Q^m) = \{\theta\} \times l_p, \quad (5.8)$$

*and*

$$E_2(Q^m) = \emptyset. \quad (5.9)$$

(iii) $m > 1$ *and m is even. Then,*

(a) *The generalized critical points of* $Q^m$ *satisfy*

$$\partial_1^C(Q^m) = \nabla_1^C(Q^m) = (\{\theta\} \times l_p) \cup \{(\bar{x},\bar{y}) \in l_p \times l_p : \bar{y} = -\bar{x}^2\}, \quad (5.10)$$

*and*

$$\partial_2^C(Q^m) = \nabla_2^C(Q^m) = \{(\bar{x},\bar{y}) \in l_p \times l_p : \bar{y} = -\bar{x}^2\}. \quad (5.11)$$

(b) *The solution sets to the ordered variational inequality problems satisfy*

$$E_1(Q^m) = (\{\theta\} \times K) \cup \{(\bar{x},\bar{y}) \in l_p \times l_p : \bar{y} = -\bar{x}^2\}. \quad (5.12)$$

*and*

$$E_2(Q^m) = \{(\bar{x},\bar{y}) \in l_p \times l_p : \bar{y} = -\bar{x}^2\}. \quad (5.13)$$

*Proof*. Recall that, for given $\bar{x} = (\bar{t}_1, \bar{t}_2, \ldots) \in l_p$ and $\bar{y} = (\bar{s}_1, \bar{s}_2, \ldots) \in l_p$, by Theorem 3.1, we have

$$\nabla_1(Q^m)(\bar{x},\bar{y}) = \begin{pmatrix} 2m\bar{t}_1(\bar{t}_1^2 + \bar{s}_1)^{m-1} & 0 & \ldots \\ 0 & 2m\bar{t}_2(\bar{t}_2^2 + \bar{s}_2)^{m-1} & 0 \\ \vdots & 0 & \ddots \end{pmatrix}$$

$$\nabla_2(Q^m)(\bar{x},\bar{y}) = \begin{pmatrix} m(\bar{t}_1^2+\bar{s}_1)^{m-1} & 0 & \dots \\ 0 & m(\bar{t}_2^2+\bar{s}_2)^{m-1} & 0 \\ \vdots & 0 & \ddots \end{pmatrix}.$$

Proof of (i). In particular, let $m = 1$. By (3.2) and (3.3) in Theorem 3.1, for given $\bar{x} = (\bar{t}_1, \bar{t}_2, \dots) \in l_p$ and $\bar{y} = (\bar{s}_1, \bar{s}_2, \dots) \in l_p$, we have

$$\nabla_1(Q)(\bar{x},\bar{y}) = \begin{pmatrix} 2\bar{t}_1 & 0 & \dots \\ 0 & 2\bar{t}_2 & 0 \\ \vdots & 0 & \ddots \end{pmatrix} \quad \text{and} \quad \nabla_2(Q)(\bar{x},\bar{y}) = \begin{pmatrix} 1 & 0 & \dots \\ 0 & 1 & 0 \\ \vdots & 0 & \ddots \end{pmatrix}.$$

And $\quad \partial_1 Q(\bar{x},\bar{y})(x) = \nabla_1 Q\ (\bar{x},\bar{y})(x) = \{2\bar{t}_n t_n\}_{n=1}^{\infty}$, for any $x = \{t_n\}_{n=1}^{\infty} \in l_p \backslash \{\theta\}$,

$$\partial_2 Q(\bar{x},\bar{y})(y) = \nabla_2 Q\ (\bar{x},\bar{y})(y) = \{s_n\}_{n=1}^{\infty}, \text{ for any } y = \{s_n\}_{n=1}^{\infty} \in l_p \backslash \{\theta\}.$$

Solving $\nabla_1(Q)(\bar{x},\bar{y}) = \theta_{XZ}$ and $\nabla_2(Q)(\bar{x},\bar{y}) = \theta_{YZ}$, this implies that $\partial_1^C(Q) = \nabla_1^C(Q) = \{\theta\} \times l_p$ and $\partial_2^C(Q) = \nabla_2^C(Q) = \emptyset$, which proves (5.4) in part (i). Meanwhile, by Proposition 4.7, the second equation in (5.5) follows from the second equation of (5.4) immediately. To completely prove part (i), we only need to prove the first equation in (5.5). By the first equation in (5.4) and by (5.3), for any $(\theta, \bar{y}) \in \{\theta\} \times l_p$, we have that

$$Q(\theta, \bar{y}) = \{\bar{s}_n\}_{n=1}^{\infty}, \text{ for any } \bar{y} = (\bar{s}_1, \bar{s}_2, \dots) \in l_p.$$

Then, for any $x = \{t_n\}_{n=1}^{\infty} \in l_p$, by definition of $Q$, we have

$$Q(x, \bar{y}) = \{t_n^2 + \bar{s}_n\}_{n=1}^{\infty} \succcurlyeq_K \{\bar{s}_n\}_{n=1}^{\infty} = Q(\theta, \bar{y}).$$

This implies that

$$(\theta, \bar{y}) \in E_1^{\wedge}(Q) \subseteq E_1(Q), \text{ for any } \bar{y} = (\bar{s}_1, \bar{s}_2, \dots) \in l_p.$$

By the first equation in (5.4) and by Theorem 3.1, this proves the first equation in (5.5). The second equation in (5.5) follows from the second equation in (5.4) immediately, which completes the proof or (i).

Proof of part (ii). By (3.3) and (3.4) in Theorem 3.1, for $\bar{x} = (\bar{t}_1, \bar{t}_2, \dots)$ and $\bar{y} = (\bar{s}_1, \bar{s}_2, \dots) \in l_p$, we have

$$\nabla_1(Q^m)(\bar{x},\bar{y}) = \begin{pmatrix} 2m\bar{t}_1(\bar{t}_1^2+\bar{s}_1)^{m-1} & 0 & \dots \\ 0 & 2m\bar{t}_2(\bar{t}_2^2+\bar{s}_2)^{m-1} & 0 \\ \vdots & 0 & \ddots \end{pmatrix},$$

$$\nabla_2(Q^m)(\bar{x},\bar{y}) = \begin{pmatrix} m(\bar{t}_1^2+\bar{s}_1)^{m-1} & 0 & \dots \\ 0 & m(\bar{t}_2^2+\bar{s}_2)^{m-1} & 0 \\ \vdots & 0 & \ddots \end{pmatrix}.$$

More precisely,

$$\partial_1 Q^m(\bar{x},\bar{y})(x) = \nabla_1(Q^m)\ (\bar{x},\bar{y})(x) = \{2m\bar{t}_n(\bar{t}_n^2+\bar{s}_n)^{m-1} t_n\}_{n=1}^{\infty}, \text{ for } x = \{t_n\}_{n=1}^{\infty} \in l_p \backslash \{\theta\}, \quad (5.14)$$

$$\partial_2 Q^m(\bar{x},\bar{y})(y) = \nabla_2(Q^m)\ (\bar{x},\bar{y})(y) = \{m(\bar{t}_n^2+\bar{s}_n)^{m-1} s_n\}_{n=1}^{\infty}, \text{ for } y = \{s_n\}_{n=1}^{\infty} \in l_p \backslash \{\theta\}. \quad (5.15)$$

Then, by $\nabla_1(Q^m)(\bar{x},\bar{y}) = \theta_{XZ}$ and $\nabla_2(Q^m)(\bar{x},\bar{y}) = \theta_{YZ}$, (5.6) and (5.7) follow from (5.14) and (5.15), respectively. Next, we prove (5.8). For any $(\theta,\bar{y}) \in \{\theta\} \times l_p$, since *m* is odd, one has

$$Q^m(\theta,\bar{y}) = \bar{y}^m \preccurlyeq_K (x^2 + \bar{y})^m = Q^m(x,\bar{y}), \text{ for any } x \in l_p.$$

This implies that

$$(\theta,\bar{y}) \in E_1^{\wedge}(Q^m) \subseteq E_1(Q^m), \text{ for any } (\theta,\bar{y}) \in \{\theta\} \times l_p.$$

This shows that

$$\{\theta\} \times l_p \subseteq E_1(Q^m). \tag{5.16}$$

For given $(\bar{x},\bar{y}) \in \partial_1^C(Q^m) = \nabla_1^C(Q^m)$ with $\bar{y} = -\bar{x}^2$ and $(\bar{x},\bar{y}) \neq (\theta,\theta)$, by $\nabla_1(Q^m)(\bar{x},\bar{y}) = \theta_{XZ}$, there is a positive integer *k* such that $\bar{s}_k \neq 0$ and $\bar{t}_k \neq 0$ satisfying $\bar{t}_k^2 + \bar{s}_k = 0$. For real number *t*, let

$$x = \left\{\bar{t}_n + \frac{t}{n}\right\}_{n=1}^{\infty} \in l_p.$$

By $\bar{y} = -\bar{x}^2$, we calculate

$$Q^m(x,\bar{y}) = \left\{\left(\left(\bar{t}_n + \frac{t}{n}\right)^2 + \bar{s}_n\right)^m\right\}_{n=1}^{\infty} = \left\{\left(\frac{2\bar{t}_n t}{n} + \left(\frac{t}{n}\right)^2\right)^m\right\}_{n=1}^{\infty}.$$

By $\bar{y} = -\bar{x}^2$ again, we have

$$Q^m(\bar{x},\bar{y}) = \{(\bar{t}_n^2 + \bar{s}_n)^m\}_{n=1}^{\infty} = \theta.$$

In particular, by $\bar{s}_k \neq 0$ and $\bar{t}_k \neq 0$, the $k^{th}$ terms of $Q^m(x,\bar{y})$ and $Q^m(\bar{x},\bar{y})$ satisfy that

$$t > \max\{0, -2k\bar{t}_k\} \text{ or } t < \min\{0, -2k\bar{t}_k\} \Longrightarrow \left(\frac{2\bar{t}_k t}{k} + \left(\frac{t}{k}\right)^2\right)^m = \left(2\bar{t}_k + \frac{t}{k}\right)^m \left(\frac{t}{k}\right)^m > 0 = \left(\bar{t}_k^2 + \bar{s}_k\right)^m;$$

$$\min\{0, -2k\bar{t}_k\} < t < \max\{0, -2k\bar{t}_k\} \Longrightarrow \left(\frac{2\bar{t}_k t}{k} + \left(\frac{t}{k}\right)^2\right)^m = \left(2\bar{t}_k + \frac{t}{k}\right)^m \left(\frac{t}{k}\right)^m < 0 = \left(\bar{t}_k^2 + \bar{s}_k\right)^m.$$

For $x = \left\{\bar{t}_n + \frac{t}{n}\right\}_{n=1}^{\infty} \in l_p$ with $\bar{y} = -\bar{x}^2$, the above two inequalities respectively imply

$$Q^m(x,\bar{y}) \text{ is } \begin{cases} \not\preccurlyeq_K Q^m(\bar{x},\bar{y}), \\ \not\succcurlyeq_K Q^m(\bar{x},\bar{y}). \end{cases}$$

Hence, for given $(\bar{x},\bar{y}) \in l_p \times l_p$ we have

$$\bar{y} = -\bar{x}^2 \text{ and } (\bar{x},\bar{y}) \neq (\theta,\theta) \implies (\bar{x},\bar{y}) \notin E_1(Q^m). \tag{5.17}$$

Then, (5.8) follows from Proposition 4.7, (5.6), (5.16) and (5.17). Then, we prove (5.9). Let $(\bar{x},\bar{y}) \in \partial_2^C(Q^m) = \nabla_2^C(Q^m)$ with $\bar{y} = -\bar{x}^2$. For real number *t*, let

$$y = \left\{\bar{s}_n + \frac{t}{n}\right\}_{n=1}^{\infty} \in l_p.$$

By $\bar{y} = -\bar{x}^2$, we calculate

$$Q^m(\bar{x}, y) = \left\{\left(\bar{t}_n^2 + \bar{s}_n + \frac{t}{n}\right)^m\right\}_{n=1}^{\infty} = \left\{\left(\frac{t}{n}\right)^m\right\}_{n=1}^{\infty}.$$

This implies that

$$t > 0 \implies Q^m(\bar{x}, y) = \left\{\left(\frac{t}{n}\right)^m\right\}_{n=1}^{\infty} \succ_K \theta = Q^m(\bar{x}, \bar{y}),$$

and
$$t < 0 \implies Q^m(\bar{x}, y) = \left\{\left(\frac{t}{n}\right)^m\right\}_{n=1}^{\infty} \prec_K \theta = Q^m(\bar{x}, \bar{y}).$$

For this arbitrarily given $(\bar{x}, \bar{y}) \in \partial_2^C(Q^m) = \nabla_2^C(Q^m)$ with $\bar{y} = -\bar{x}^2$, we obtain that

$$(\bar{x}, \bar{y}) \in \partial_2^C(Q^m) \quad \text{and} \quad (\bar{x}, \bar{y}) \notin E_2(Q^m). \tag{5.18}$$

Then, (5.9) follows from Proposition 4.7, (5.7) and (5.18).

Proof of (iii). (5.10) and (5.11) follow from (5.14) and (5.15), respectively. Next, we prove (5.12). For any $(\theta, \bar{y}) \in \{\theta\} \times K$, since $m$ is even, for any given $\bar{y} \in K$, one has

$$Q^m(\theta, \bar{y}) = \bar{y}^m \preccurlyeq_K (x^2 + \bar{y})^m = Q^m(x, \bar{y}), \text{ for any } x \in l_p.$$

This implies that

$$(\theta, \bar{y}) \in E_1^{\wedge}(Q^m) \subseteq E_1(Q^m), \text{ for any } (\theta, \bar{y}) \in \{\theta\} \times K. \tag{5.19}$$

Let $(\theta, \bar{y}) \in \{\theta\} \times l_p$ with $(\theta, \bar{y}) \notin \{\theta\} \times K$. Then there is a positive integer $k$ such that $\bar{s}_k < 0$. Let $x = \{t_n\}_{n=1}^{\infty} \in l_p$ with $t_n = 0$, for $n \neq k$ and $t_k \neq 0$. Since $m$ is even, we have

$$\left(t_k^2 + \bar{s}_k\right)^m \text{ is } \begin{cases} < (0 + \bar{s}_k)^m, & \text{if } 0 < |t_k| < \sqrt{|\bar{s}_k|}, \\ > (0 + \bar{s}_k)^m, & \text{if } |t_k| > \sqrt{2|\bar{s}_k|}. \end{cases}$$

This implies that, for such a $x = \{t_n\}_{n=1}^{\infty} \in l_p$, for $(\theta, \bar{y}) \notin \{\theta\} \times K$, we have

$$Q^m(x, \bar{y}) \text{ is } \begin{cases} \not\succcurlyeq_K Q^m(\theta_p, \bar{y}), & \text{if } 0 < |t_k| < \sqrt{|\bar{s}_k|}, \\ \not\preccurlyeq_K Q^m(\theta_p, \bar{y}), & \text{if } |t_k| > \sqrt{2|\bar{s}_k|}. \end{cases}$$

Hence, we obtain that

$$(\theta, \bar{y}) \notin E_1(Q^m), \text{ for any } (\theta, \bar{y}) \notin \{\theta\} \times K. \tag{5.20}$$

Let $(\bar{x}, \bar{y}) \in (l_p \times l_p)$ with $\bar{x}^2 + \bar{y} = \theta$. Then, for any $x = \{t_n\}_{n=1}^{\infty} \in l_p$, we have (since $m$ is even)

$$(t_n^2 + \bar{s}_n)^m \geq 0 = (\bar{t}_n^2 + \bar{s}_n)^m, \text{ for every } n = 1, 2, \dots .$$

This implies that, for any $(\bar{x}, \bar{y}) \in (l_p \times l_p)$ with $\bar{x}^2 + \bar{y} = \theta$, we have

$$Q^m(x, \bar{y}) \succcurlyeq_K Q^m(\bar{x}, \bar{y}), \text{ for any } x = \{t_n\}_{n=1}^{\infty} \in l_p. \tag{5.21}$$

By definition, (5.21) implies that

$$(\bar{x}, \bar{y}) \in E_1^\wedge(Q^m) \subseteq E_1(Q^m), \text{ for any } (\bar{x}, \bar{y}) \in (l_p \times l_p) \text{ with } \bar{x}^2 + \bar{y} = \theta. \tag{5.22}$$

By Proposition 4.7, (5.12) is proved by following (5.10), (5.19), (5.20) and (5.22).

To complete the proof of part (iii), finally, we prove (5.13). As a matter of fact, similar to the proof of (5.21), for $(\bar{x}, \bar{y}) \in (l_p \times l_p)$ with $\bar{x}^2 + \bar{y} = \theta$, we have

$$Q^m(\bar{x}, y) \succcurlyeq_K Q^m(\bar{x}, \bar{y}), \text{ for any } y = \{s_n\}_{n=1}^\infty \in l_p.$$

This implies that

$$(\bar{x}, \bar{y}) \in E_2^\wedge(Q^m) \subseteq E_2(Q^m), \text{ for any } (\bar{x}, \bar{y}) \in (l_p \times l_p) \text{ with } \bar{x}^2 + \bar{y} = \theta. \tag{5.23}$$

By Proposition 4.7 and (5.11), (5.13) is proved by (5.23). □

## 6. Some Duality Mappings in $l_p$

In this section, we consider some duality mappings in the uniformly convex and uniformly smooth Banach space $l_p$, in which $1 < p, q < \infty$ with $\frac{1}{p} + \frac{1}{q} = 1$, We study their partial Gâteaux differentiability in $l_p$, find their generalized partially critical points, and calculate their partial ordered extrema. To start this section, we recall some notations and concepts used in [21].

For $1 < p < \infty$, let $(l_p, \|\cdot\|_p, \preccurlyeq_K)$ be the partially ordered Banach spaces, in which the partial order $\preccurlyeq_K$ is generated by the positive cone $K$ in $l_p$, and let $(l_q, \|\cdot\|_q, \preccurlyeq_M)$ be the partially ordered Banach spaces, in which the partial order $\preccurlyeq_M$ is generated by the positive cone $M$ in $l_q$.

Let $\mathbb{S}$ denote the linear space of real sequences with term-by-term addition and scalar multiplication. Let $x = (t_1, t_2, \ldots) \in \mathbb{S}$. Let arbitrarily given $\bar{x} = (\bar{t}_1, \bar{t}_2, \ldots) \in \mathbb{S}$. Define an operator $R(\bar{x})(\cdot): \mathbb{S} \to \mathbb{S}$, for $v = \{s_n\}_{n=1}^\infty \in \mathbb{S}$ by $R(\bar{x})(v) := \{\lambda_n\}_{n=1}^\infty$, in which, for each $n = 1, 2, \ldots$, $\lambda_n \in \mathbb{R}$ and

$$\lambda_n = \begin{cases} 0, & \text{if } \bar{t}_n = 0, \\ \frac{s_n}{\bar{t}_n}, & \text{if } \bar{t}_n \neq 0. \end{cases}$$

The sequence $R(\bar{x})(v) = \{\lambda_n\}_{n=1}^\infty \in \mathbb{S}$ is called the ratio sequence of $v$ with respect to $\bar{x}$.

Let $f$: $l_p \to l_q$ be a single-valued mapping. If there is a continuous and strictly increasing function $\varphi: [0, \infty) \to [0, \infty)$ satisfying $\varphi(0) = 0$ and $\varphi(t) \to \infty$, as $t \to \infty$, such that

$$\|f(x)\|_q = \varphi(\|x\|_p), \text{ for any } x \in l_p, \tag{6.1}$$

then, $f$ is called a duality mapping in $l_p$ with gauge function $\varphi$. In case, if the sign "=" in (6.1) is replaced by "≤", then $f$ is called a generalized duality mapping in $l_p$ with the gauge function $\varphi$. Especially, the normalized duality mapping $J$: $l_p \to l_q$ is a special duality mapping in $l_p$ with gauge function $\varphi$ defined by $\varphi(t) = t$, for any $t \geq 0$. More precisely, the normalized duality mapping $J$: $l_p \to l_q$ is a continuous single-valued mapping, which satisfies that $J(\theta) = \theta$ and for any $x = (t_1, t_2, \ldots) \in l_p \backslash \{\theta\}$,

$$J(x) = \left(\frac{|t_1|^{p-1}\text{sign}(t_1)}{\|x\|_p^{p-2}}, \frac{|t_2|^{p-1}\text{sign}(t_2)}{\|x\|_p^{p-2}}, \ldots\right). \tag{6.2}$$

In sections 7 of [21], a generalized duality mapping $\mathfrak{J}: l_p \to l_q$ is defined:

$$\mathfrak{J}(y) = \left\{\frac{|s_n|^{p-1}\text{sign}(s_n)}{1+|s_n|}\right\}_{n=1}^{\infty}, \text{ for any } y = (s_1, s_2, \ldots) \in l_p.$$

We review some results about the Gâteaux differentiability of the mappings $J$ and $\mathfrak{J}$ in $l_p$ proved in [21], which will be used in this section to study the partial Gâteaux differentiability of the sum of $J$ and $\mathfrak{J}$ in $l_p$.

**Partial of Theorem 3.2 in [21]**. *The normalized duality mapping $J$: $l_p \to l_q$ has differentiability below*:

(i) *For $1 < p < \infty$, $J$ is Gâteaux differentiable at $\theta$ and*

$$J'(\theta)(v) = J(v), \textit{ for any } v \in l_p\backslash\{\theta\}. \tag{6.3}$$

(ii) *For $2 < p < \infty$, $J$ is Gâteaux differentiable on $l_p\backslash\{\theta\}$. More precisely, let $\bar{x} = (\bar{t}_1, \bar{t}_2, \ldots) \in l_p\backslash\{\theta\}$. For any $v = \{s_n\}_{n=1}^{\infty} \in l_p\backslash\{\theta\}$ with ratio sequence $\{\lambda_n\}_{n=1}^{\infty}$ with respect to $\bar{x}$, the Gâteaux derivative of $J$ at $\bar{x}$ in direction $v$ is*

$$J'(\bar{x})(v) = J(\bar{x})\left((p-1)\begin{pmatrix}\lambda_1 & 0 & \ldots \\ 0 & \lambda_2 & 0 \\ \vdots & 0 & \ddots\end{pmatrix} - \frac{(p-2)\sum_{i=1}^{\infty}\lambda_i|\bar{t}_i|^p}{\sum_{i=1}^{\infty}|\bar{t}_i|^p}\begin{pmatrix}1 & 0 & \ldots \\ 0 & 1 & 0 \\ \vdots & 0 & \ddots\end{pmatrix}\right)$$

$$= \left\{\frac{|\bar{t}_n|^{p-1}\text{sign}(\bar{t}_n)}{\|\bar{x}\|_p^{p-2}}\left((p-1)\lambda_n - \frac{(p-2)\sum_{i=1}^{\infty}\lambda_i|\bar{t}_i|^p}{\|\bar{x}\|_p^p}\right)\right\}_{n=1}^{\infty}. \tag{6.4}$$

**Corollary 7.3 in [21]**. *Let $3 < p < \infty$. $\mathfrak{J}$ is Gâteaux differentiable in $l_p$ and, for any $\bar{y} = (\bar{s}_1, \bar{s}_2, \ldots) \in l_p$, the Gâteaux derivative $\mathfrak{J}'(\bar{y})$ satisfies that for any $u = \{s_n\}_{n=1}^{\infty} \in l_p\backslash\{\theta\}$,*

$$\mathfrak{J}'(\bar{y})(u) = \left\{\frac{(p-1)|\bar{s}_n|^{p-2}+(p-2)|\bar{s}_n|^{p-1}}{(1+|\bar{s}_n|)^2}s_n\right\}_{n=1}^{\infty} \in l_q. \tag{6.5}$$

By using the normalized duality mapping $J$ and the generalized duality mapping $\mathfrak{J}$, we define a two-variable mapping $F$: $l_p \times l_p \to l_q$ by

$$F(u, v) = J(u) + \mathfrak{J}(v), \text{ for any } v = \{t_n\}_{n=1}^{\infty}, u = \{s_n\}_{n=1}^{\infty} \in l_p. \tag{6.6}$$

**Theorem 6.1**. *Let $3 < p < \infty$ with $\frac{1}{p} + \frac{1}{q} = 1$. Let $\bar{x} = (\bar{t}_1, \bar{t}_2, \ldots) \in l_p$ and $\bar{y} = (\bar{s}_1, \bar{s}_2, \ldots) \in l_p$. The mapping $F$ defined in* (6.6) *has the following partial derivatives at point $(\bar{x}, \bar{y})$: for any $v = \{t_n\}_{n=1}^{\infty} \in l_p\backslash\{\theta\}$ with ratio sequence $\{\lambda_n\}_{n=1}^{\infty}$ with respect to $\bar{x}$ and for any $u = \{s_n\}_{n=1}^{\infty} \in l_p\backslash\{\theta\}$,*

(i) $\partial_1 F(\theta, \bar{y})(v) = J(v)$;

(ii) $$\partial_1 F(\bar{x}, \bar{y})(v) = \left\{\frac{|\bar{t}_n|^{p-1}\text{sign}(\bar{t}_n)}{\|\bar{x}\|_p^{p-2}}\left((p-1)\lambda_n - \frac{(p-2)\sum_{i=1}^{\infty}\lambda_i|\bar{t}_i|^p}{\|\bar{x}\|_p^p}\right)\right\}_{n=1}^{\infty}, \textit{ if } \bar{x} \neq \theta;$$

(iii) $\partial_2 F(\bar{x},\bar{y})(u) = \left\{\frac{(p-1)|\bar{s}_n|^{p-2}+(p-2)|\bar{s}_n|^{p-1}}{(1+|\bar{s}_n|)^2} s_n\right\}_{n=1}^{\infty}$.

*Proof*. By Lemma 2.10, for any $v = \{t_n\}_{n=1}^{\infty} \in l_p\backslash\{\theta\}$ with ratio sequence $\{\lambda_n\}_{n=1}^{\infty}$ with respect to $\bar{x}$ and for any $u = \{s_n\}_{n=1}^{\infty} \in l_p\backslash\{\theta\}$, by Theorem 3.2 and Corollary 7.3 in [21], we have

$$\partial_1 F(\theta,\bar{y})(v) = J'(\theta)(v) = J(v);$$

$$\partial_1 F(\bar{x},\bar{y})(v) = J'(\bar{x})(v) = \left\{\frac{|\bar{t}_n|^{p-1}\mathrm{sign}(\bar{t}_n)}{\|\bar{x}\|_p^{p-2}}\left((p-1)\lambda_n - \frac{(p-2)\sum_{i=1}^{\infty}\lambda_i|\bar{t}_i|^p}{\|\bar{x}\|_p^p}\right)\right\}_{n=1}^{\infty}, \text{ if } \bar{x} \neq \theta;$$

$$\partial_2 F(\bar{x},\bar{y})(u) = \mathfrak{J}'(\bar{y})(u) = \left\{\frac{(p-1)|\bar{s}_n|^{p-2}+(p-2)|\bar{s}_n|^{p-1}}{(1+|\bar{s}_n|)^2} s_n\right\}_{n=1}^{\infty}. \qquad \square$$

**Theorem 6.2**. *Let* $3 < p < \infty$. *The generalized partially critical points of the mapping* $F$ *defined in* (6.6) *have the following properties.*

(i) $\partial_1^C(F) = \emptyset$ *and* $E_1(F) = \emptyset$;
(ii) $\partial_2^C(F) = l_p \times \{\theta\}$ *and* $E_2(F) = \emptyset$;

*Proof*. Proof of (i). By (i) of Theorem 6.1, we immediately have that

$$(\theta,\bar{y}) \notin \partial_1^C(F), \text{ for any } \bar{y} \in l_p. \tag{6.7}$$

Let arbitrarily given $(\bar{x},\bar{y}) \in l_p \times l_p$ with $\bar{x} \neq \theta$ and $\bar{x} = (\bar{t}_1, \bar{t}_2, \dots)$. There is a positive number $k$ with $\bar{t}_k \neq 0$. Let $v = \{t_n\}_{n=1}^{\infty} \in l_p\backslash\{\theta\}$ with ratio sequence $\{\lambda_n\}_{n=1}^{\infty}$ with respect to $\bar{x}$ satisfying $t_n = 0$, for all $n \neq k$ and $t_k = \bar{t}_k$. This implies that $\lambda_k = 1$ and $\lambda_n = 0$, for all $n \neq k$. By part (ii) of Theorem 6.1, The $k$th entry of $\frac{\partial F(x,y)}{\partial x}(\bar{x},\bar{y})(v)$ is

$$\frac{|\bar{t}_k|^{p-1}\mathrm{sign}(\bar{t}_k)}{\|\bar{x}\|_p^{p-2}}\left((p-1) - \frac{(p-2)|\bar{t}_k|^p}{\|\bar{x}\|_p^p}\right).$$

Notice that $\|\bar{x}\|_p^p = \sum_{i=1}^{\infty}|\bar{t}_i|^p$. By $3 < p < \infty$, we estimate

$$\frac{(p-2)|\bar{t}_k|^p}{\|\bar{x}\|_p^p} \leq p-2 < p-1.$$

This implies that

$$\partial_1 F(\bar{x},\bar{y})(v) \neq \theta, \text{ for any } (\bar{x},\bar{y}) \in l_p \times l_p \text{ with } \bar{x} \neq \theta. \tag{6.8}$$

By (6.7) and (6.8), we obtain that

$$(\bar{x},\bar{y}) \notin \partial_1^C(F), \text{ for any } (\bar{x},\bar{y}) \in l_p \times l_p \text{ with } \bar{x} \neq \theta.$$

This proves the first equation of (i) of this theorem. By Proposition 4.7, this implies the second equation of (i). Then, we prove (ii) of this theorem. By (ii) of Theorem 6.1, it is clear that

$$l_p \times \{\theta\} \subseteq \partial_2^C(F). \tag{6.9}$$

For any $(\bar{x}, \bar{y}) \in l_p \times l_p$ with $\bar{y} \neq \theta$ and $\bar{y} = (\bar{s}_1, \bar{s}_2, \dots) \in l_p$. There is a positive number $k$ with $\bar{s}_k \neq 0$. Let $u = \{s_n\}_{n=1}^{\infty} \in l_p \backslash \{\theta\}$ with $s_k = \bar{s}_k$. By part (iii) of Theorem 6.1, The $k$th entry of $\partial_2 F(\bar{x}, \bar{y})(u)$ is

$$\frac{(p-1)|\bar{s}_k|^{p-2} + (p-2)|\bar{s}_k|^{p-1}}{(1+|\bar{s}_k|)^2} \bar{s}_k \neq 0.$$

This implies that, for the given $(\bar{x}, \bar{y}) \in l_p \times l_p$ with $\bar{y} \neq \theta$, we have

$$\partial_2 F(\bar{x}, \bar{y})(u) \neq \theta, \text{ for } u = \{s_n\}_{n=1}^{\infty} \in l_p \backslash \{\theta\} \text{ with } s_k = \bar{s}_k.$$

Hence, we obtain that

$$(\bar{x}, \bar{y}) \notin \partial_2^C(F), \text{ for any } (\bar{x}, \bar{y}) \in l_p \times l_p \text{ with } \bar{y} \neq \theta. \tag{6.10}$$

The first equation of part (ii) of this theorem is proved by (6.9) and (6.10). For any given $(\bar{x}, \theta) \in l_p \times \{\theta\}$ with $\bar{x} = (\bar{t}_1, \bar{t}_2, \dots) \neq \theta$, we have

$$F(\bar{x}, \theta) = J(\bar{x}) + \mathfrak{J}(\theta) = J(\bar{x}) = \left\{\frac{|\bar{t}_n|^{p-1}\text{sign}(\bar{t}_n)}{\|\bar{x}\|_p^{p-2}}\right\}_{n=1}^{\infty}. \tag{6.11}$$

For this fixed $\bar{x} \neq \theta$, for $u = \{s_n\}_{n=1}^{\infty} \in l_p$, we have

$$F(\bar{x}, u) = J(\bar{x}) + \mathfrak{J}(u)$$

$$= \left\{\frac{|\bar{t}_n|^{p-1}\text{sign}(\bar{t}_n)}{\|\bar{x}\|_p^{p-2}}\right\}_{n=1}^{\infty} + \left\{\frac{|s_n|^{p-1}\text{sign}(s_n)}{1+|s_n|}\right\}_{n=1}^{\infty}$$

$$= \left\{\frac{|\bar{t}_n|^{p-1}\text{sign}(\bar{t}_n)}{\|\bar{x}\|_p^{p-2}} + \frac{|s_n|^{p-1}\text{sign}(s_n)}{1+|s_n|}\right\}_{n=1}^{\infty}.$$

For any given positive number $n$, we have

$$\frac{|\bar{t}_n|^{p-1}\text{sign}(\bar{t}_n)}{\|\bar{x}\|_p^{p-2}} + \frac{|s_n|^{p-1}\text{sign}(s_n)}{1+|s_n|} \quad \text{is} \quad \begin{cases} > \frac{|\bar{t}_n|^{p-1}\text{sign}(\bar{t}_n)}{\|\bar{x}\|_p^{p-2}}, & \text{if } s_n > 0, \\ < \frac{|\bar{t}_n|^{p-1}\text{sign}(\bar{t}_n)}{\|\bar{x}\|_p^{p-2}}, & \text{if } s_n < 0. \end{cases}$$

This implies that

$$F(\bar{x}, u) \text{ is} \begin{cases} \ngeq_M F(\bar{x}, \theta), & \text{if } s_n < 0, \\ \nleq_M F(\bar{x}, \theta), & s_n > 0. \end{cases} \tag{6.12}$$

Then, for $\bar{x} \neq \theta$, by comparing (6.13) and (6.11), this implies that $F$ does not have neither partial ordered-minimum value, nor partial ordered-maximum value at point $(\bar{x}, \theta)$ with respect to the second variable.

Let $\bar{x} = \theta$. For $u = \{s_n\}_{n=1}^{\infty} \in l_p$, we have

$$F(\theta, u) = J(\theta) + \mathfrak{J}(u)$$

$$= \{0\}_{n=1}^{\infty} + \left\{\frac{|s_n|^{p-1}\text{sign}(s_n)}{1+|s_n|}\right\}_{n=1}^{\infty}$$

$$= \left\{\frac{|s_n|^{p-1}\text{sign}(s_n)}{1+|s_n|}\right\}_{n=1}^{\infty}.$$

For any given positive number $n$, we have

$$\frac{|\bar{t}_n|^{p-1}\text{sign}(\bar{t}_n)}{\|\bar{x}\|_p^{p-2}} + \frac{|s_n|^{p-1}\text{sign}(s_n)}{1+|s_n|} \quad \text{is} \begin{cases} > 0, & \text{if } s_n > 0, \\ < 0, & \text{if } s_n < 0. \end{cases}$$

This implies that

$$F(\theta, u) \text{ is} \begin{cases} \not\succcurlyeq_M F(\theta,\theta), & \text{if } s_n < 0, \\ \not\preccurlyeq_M F(\theta,\theta), & s_n > 0. \end{cases} \tag{6.13}$$

Then, for $\bar{x} = \theta$, by comparing (6.13) and (6.12), this implies that $F$ does not have neither partial ordered-minimum value, nor partial ordered-maximum value at point $(\theta, \theta)$ with respect to the second variable, which proves $E_2(F) = \emptyset$. □

## 7. Ordered Variational Inequality Problems

In this section, we show that the ordered optimization problems studied in section 4 are equivalently converted to ordered variational inequality problems in partially ordered Banach spaces. We review the concepts of ordered variational inequality problems, as usual.

Following the notations in sections 2 and 4, let $(X, \|\cdot\|_X)$ and $(Y, \|\cdot\|_Y)$ be Banach spaces and $(Z, \|\cdot\|_Z, \preccurlyeq_M)$ a partially ordered Banach spaces, in which the partial order $\preccurlyeq_M$ is generated by a nonempty closed convex and pointed cone $M$ in $Z$. Let $\theta_X$, $\theta_Y$ and $\theta_Z$ be the origins of $X$, $Y$ and $Z$ respectively. Let $A$ and $B$ be nonempty convex open subsets in $X$ and $Y$ respectively.

**Definition 7.1**. Let $F: A \times B \to Z$ be a single-valued mapping. An ordered variational inequality problem associated with $F$, $A$ and $B$ with respect to the first variable, denoted by $\text{OVI}_1(F; A, B)$, is to find a point $(\bar{x}, \bar{y}) \in A \times B$ such that

$$F(\bar{x}, \bar{y}) \preccurlyeq_M F(x, \bar{y}), \text{ for every } x \in A. \tag{7.1}$$

A point $(\bar{x}, \bar{y}) \in A \times B$ satisfying (7.1) is called a solution to the ordered variational inequality problem $\text{OVI}_1(F; A, B)$. The set of all solutions to the problem $\text{OVI}_1(F; A, B)$ is denoted by $\text{SOVI}_1(F; A, B)$.

The ordered variational inequality problem with respect to the second variable $\text{OVI}_2(F; A, B)$ is similarly defined to find a point $(\bar{x}, \bar{y}) \in A \times B$ such that

$$F(\bar{x}, \bar{y}) \preccurlyeq_M F(\bar{x}, y), \text{ for every } y \in B. \tag{7.2}$$

The set of all solutions to the problem $\text{OVI}_2(F; A, B)$ is denoted by $\text{SOVI}_2(F; A, B)$.

**Remarks 7.2.** Comparing (7.1) with (4.1) and (7.2) with (4.2), we have

(i) $\text{SOVI}_1(F; A, B) = E_1^{\wedge}(F, A, B)$;
(ii) $\text{SOVI}_2(F; A, B) = E_2^{\wedge}(F, A, B)$.

Then, Proposition 4.5 is rewritten with respect to ordered variational inequality problems.

**Proposition 7.3**. *Let $F: A \times B \to Z$ be a single-valued mapping.*

(i) *Suppose that F is partial Gâteaux differentiable on A× B with respect to every variable. Then*

$$\text{SOVI}_1(F; A, B) \subseteq \partial_1^C F \quad \textit{and} \quad \text{SOVI}_2(F; A, B) \subseteq \partial_2^C F.$$

(ii) *Suppose that F is partial Fréchet differentiable on A× B with respect to every variable. Then*

$$\text{SOVI}_2(F; A, B) \subseteq \nabla_1^C F \quad \textit{and} \quad \text{SOVI}_2(F; A, B) \subseteq \nabla_2^C F.$$

## 8. Ordered Monotonicity in Partially Ordered Banach Spaces

let $(X, \|\cdot\|_X, \preccurlyeq_C)$, $(Y, \|\cdot\|_Y, \preccurlyeq_K)$ and $(Z, \|\cdot\|_Z, \preccurlyeq_M)$ be partially ordered Banach spaces, in which the partial orders $\preccurlyeq_C$, $\preccurlyeq_K$ and $\preccurlyeq_M$ are generated by nonempty closed convex and pointed cones $C$, $K$ and $M$ in $X$, $Y$ and $Z$ respectively. Let $A$ and $B$ be nonempty convex open subsets in $X$ and $Y$ respectively.

**Definition 8.1**. Let $F: A\times B \to Z$ be a single-valued mapping. Let $\bar{x} \in A$ and $\bar{y} \in B$. If, for any $x, v \in A$,

$$x \preccurlyeq_C v \quad \Longrightarrow \quad F(x, \bar{y}) \preccurlyeq_M F(v, \bar{y}), \tag{8.1}$$

then, $F$ is said to be order-increasing in $A$ with respect to the first variable at point $\bar{y} \in B$. The order decreasing of $F$ in $A$ with respect to the first variable at point $\bar{y} \in B$ can be similarly defined. If, for any $y, u \in A$,

$$y \preccurlyeq_K u \quad \Longrightarrow \quad F(\bar{x}, y) \preccurlyeq_M F(\bar{x}, u), \tag{8.2}$$

then, $F$ is said to be order-increasing in $B$ with respect to the second variable at point $\bar{x} \in A$. The order decreasing of $F$ in $B$ with respect to the first variable at point $\bar{x} \in A$ can be similarly defined.

**Theorem 8.2**. *Let $F: A\times B \to Z$ be a single-valued mapping. Suppose that F is Gâteaux differentiable with each variable in $A\times B$. If F is order increasing (decreasing) in A with respect to the first variable at some point $\bar{y} \in B$, then, for any given $w \in A$,*

$$\partial_1 F(w, \bar{y})(v) \succcurlyeq_M (\preccurlyeq_M)\ \theta_Z, \textit{ for any } v \in C\backslash\{\theta_X\}. \tag{8.3}$$

*In particular, if F is Fréchet differentiable with each variable in $A\times B$ and F is order increasing (decreasing) in A with respect to the first variable at some point $\bar{y} \in B$, then for any given $w \in A$,*

$$\nabla_1 F(w, \bar{y})(v) \succcurlyeq_M (\preccurlyeq_M)\ \theta_Z, \textit{ for any } v \in C\backslash\{\theta_X\}. \tag{8.4}$$

*Proof*. We first prove (8.3). Suppose that $F$ is order increasing in $A$ with respect to the first variable at some point $\bar{y} \in Y$. Since $A$ is a nonempty convex open subset in $X$, for the given $w \in A$ and $v \in C\backslash\{\theta_X\}$, there is $\delta_0 > 0$ such that, for real number $t$, we have

$$|t| < \delta_0 \quad \Longrightarrow \quad w + tv \in A.$$

By $v \in C\backslash\{\theta_X\}$, we have that

$$t > 0 \;\Longrightarrow\; w + tv \succcurlyeq_C w \quad \text{and} \quad t < 0 \;\Longrightarrow\; w + tv \preccurlyeq_C w. \tag{8.5}$$

By the condition that $F$ is order increasing in $A$ with respect to the first variable at the given point $\bar{y} \in Y$, by (8.5), we always have

$$\frac{F(w+tv,\bar{y}) - F(w,\bar{y})}{t} \succcurlyeq_M \theta_Z, \text{ for any } t \text{ with } |t| < \delta_0. \tag{8.6}$$

(8.6) is equivalent to

$$\frac{F(w+tv,\bar{y})-F(w,\bar{y})}{t} \in M, \text{ for any } t \text{ with } |t| < \delta_0. \tag{8.7}$$

Since $M$ is a nonempty closed convex and pointed cone in $Z$, (8.7) implies that

$$\partial_1 F(w,\bar{y})(v) = \lim_{t\to 0} \frac{F(w+tv,\bar{y})-F(w,\bar{y})}{t} \in M.$$

This proves (8.3). (8.4) is proved by (8.3) and the connection between partial Gâteaux and Fréchet derivatives. □

Similarly to Theorem 8.2, we consider the ordered monotonicity of $F$ with respect to the second variable.

**Corollary 8.3**. *Let $F$: $A\times B \to Z$ be a single-valued mapping. Suppose that $F$ is Gâteaux differentiable with each variable in $A\times B$. If $F$ is order increasing* (*decreasing*) *in $B$ with respect to the second variable at some point $\bar{x} \in A$, then, for any given $w \in B$,*

$$\partial_2 F(\bar{x},w)(u) \succcurlyeq_M (\preccurlyeq_M)\ \theta_Z, \text{ for any } u \in K\backslash\{\theta_Y\}. \tag{8.8}$$

*In particular, if $F$ is Fréchet differentiable with each variable in $A\times B$ and $F$ is order increasing* (*decreasing*) *in $B$ with respect to the second variable at some point $\bar{x} \in A$, then for any given $w \in B$,*

$$\nabla_2 F(\bar{x},w)(u) \succcurlyeq_M (\preccurlyeq_M)\ \theta_Z, \text{ for any } u \in K\backslash\{\theta_Y\}. \tag{8.9}$$

In particular, similarly to section 6, for $1 < p < \infty$, let $(l_p, \|\cdot\|_p, \preccurlyeq_K)$ be the partially ordered Banach spaces, in which the partial order $\preccurlyeq_K$ is generated by the positive cone $K$ in $l_p$, and let $(l_q, \|\cdot\|_q, \preccurlyeq_M)$ be the partially ordered Banach spaces, in which the partial order $\preccurlyeq_M$ is generated by the positive cone $M$ in $l_q$. In section 6 in [21], a duality mapping $\mathcal{J}$: $l_p \to l_q$ is defined, for any $x = (t_1, t_2, \dots) \in l_p$, by

$$\mathcal{J}(x) = (|t_1|^{p-1}\text{sign}(t_1),\ |t_2|^{p-1}\text{sign}(t_2), \dots). \tag{8.10}$$

One can show that for each $n$, the function $|t_n|^{p-1}\text{sign}(t_n)$ is a strictly increasing function of $t_n$ on $(-\infty, +\infty)$. This duality mapping $\mathcal{J}$ has the following differentiational properties proved in [21].

**Theorem 6.1 in [21]**. *Let* $3 < p < \infty$. *$\mathcal{J}$ is Fréchet differentiable on $l_p$ and, for any $\bar{x} = (\bar{t}_1, \bar{t}_2, \dots) \in l_p$,*

$$\nabla\mathcal{J}(\bar{x}) = (p-1)\begin{pmatrix} |\bar{t}_1|^{p-2} & 0 & \dots \\ 0 & |\bar{t}_2|^{p-2} & 0 \\ \vdots & 0 & \ddots \end{pmatrix}$$

*and* $$\nabla\mathcal{J}(\bar{x})(v) = \{t_n\}_{n=1}^{\infty}(p-1)\begin{pmatrix} |\bar{t}_1|^{p-2} & 0 & \dots \\ 0 & |\bar{t}_2|^{p-2} & 0 \\ \vdots & 0 & \ddots \end{pmatrix} = \{(p-1)|\bar{t}_n|^{p-2}t_n\}_{n=1}^{\infty} \in l_q. \tag{8.11}$$

Recall that the generalized duality mapping $\mathfrak{J}$: $l_p \to l_q$ is used in the previous sections.

$$\mathfrak{J}(y) = \left\{\frac{|s_n|^{p-1}\text{sign}(s_n)}{1+|s_n|}\right\}_{n=1}^{\infty}, \text{ for any } y = (s_1, s_2, \dots) \in l_p. \tag{8.12}$$

One can show that for each $n$, the function $\frac{|s_n|^{p-1}\text{sign}(s_n)}{1+|s_n|}$ is a strictly increasing function of $s_n$ on $(-\infty, +\infty)$.

**Theorem 7.2 in [21]**. *Let* 3 < *p* < ∞. $\mathfrak{J}$ *is Fréchet differentiable in* $l_p$ *and, for any* $\bar{x} = (\bar{t}_1, \bar{t}_2, \dots) \in l_p$,

$$\nabla\mathfrak{J}(\bar{x}) = \begin{pmatrix} \frac{\left((p-1)|\bar{t}_1|^{p-2}+(p-2)|\bar{t}_1|^{p-1}\right)}{(1+|\bar{t}_1|)^2} & 0 & \dots \\ 0 & \frac{\left((p-1)|\bar{t}_2|^{p-2}+(p-2)|\bar{t}_2|^{p-1}\right)}{(1+|\bar{t}_2|)^2} & 0 \\ \vdots & 0 & \ddots \end{pmatrix}.$$

*More precisely,* $\nabla\mathfrak{J}(\bar{x})$: $l_p \to l_q$ *is a pointwise multiplication operator in* $l_p$ *and, for* $v = \{t_n\}_{n=1}^{\infty} \in l_p$,

$$\nabla\mathfrak{J}(\bar{x})(v) = \left\{\frac{\left((p-1)|\bar{t}_n|^{p-2}+(p-2)|\bar{t}_n|^{p-1}\right)}{(1+|\bar{t}_n|)^2} t_n\right\}_{n=1}^{\infty} \in l_q, \tag{8.13}$$

Let 3 < *p* < ∞. We define a single-valued mapping $F$: $l_p \times l_p \to l_q$ as follows.

$$F(x, y) = \mathcal{J}(x) + \mathfrak{J}(y), \text{ for any } x, y \in l_p. \tag{8.14}$$

Next, we give a counter example to show that the conditions (8.4) and (8.5) are necessary conditions for *F* to be order increasing with respect to the first and the second variable, respectively.

**Example 8.5**. Let 3 < *p* < ∞. Let $F$: $l_p \times l_p \to l_q$ be defined by (8.14). For any $\bar{x} = (\bar{t}_1, \bar{t}_2, \dots) \in l_p$ and $\bar{y} = (\bar{s}_1, \bar{s}_2, \dots) \in l_p$, by (8.10) and (8.12), we have that *F* is order increasing with respect to the first variable at given point $\bar{y} \in l_p$, and *F* is order increasing with respect to the second variable at $\bar{x} \in l_p$. By (8.11) and (8.12), we have

$$\nabla_1 F(\bar{x}, \bar{y}) = \nabla\mathcal{J}(\bar{x}) = (p-1)\begin{pmatrix} |\bar{t}_1|^{p-2} & 0 & \dots \\ 0 & |\bar{t}_2|^{p-2} & 0 \\ \vdots & 0 & \ddots \end{pmatrix},$$

and

$$\nabla_2 F(\bar{x}, \bar{y}) = \nabla\mathfrak{J}(\bar{x}) = \begin{pmatrix} \frac{\left((p-1)|\bar{t}_1|^{p-2}+(p-2)|\bar{t}_1|^{p-1}\right)}{(1+|\bar{t}_1|)^2} & 0 & \dots \\ 0 & \frac{\left((p-1)|\bar{t}_2|^{p-2}+(p-2)|\bar{t}_2|^{p-1}\right)}{(1+|\bar{t}_2|)^2} & 0 \\ \vdots & 0 & \ddots \end{pmatrix}.$$

This immediately implies that

$$\nabla_1 F(\bar{x}, \bar{y})(v) \succcurlyeq_M \theta, \text{ for any } v \in K\backslash\{\theta\},$$

and

$$\nabla_2 F(\bar{x}, \bar{y})(u) \succcurlyeq_M \theta, \text{ for any } u \in K\backslash\{\theta\}.$$

Hence, conditions (8.4) and (8.9) are satisfied for *F* defined in (8.15), which verifies the results of Theorem 8.2 and Corollary 8.3.

Now, we provide other examples below to verify the results of Theorem 8.2 and Corollary 8.3 again. Let *m* be a positive integer. We define the *m*th power operator $P^m$: $l_p \to l_p$ as follows.

$$P^m(x) = \{t_n^m\}_{n=1}^{\infty}, \text{ for any } x = \{t_n\}_{n=1}^{\infty} \in l_p. \tag{8.15}$$

**Theorem 3.1 and Corollary 3.2 in [20]**. *Let m be a positive integer. Then, the m*th *power operator* $P^m$ *defined in* (8.15) *has the following differentiation properties.*

(i) $P^m$ *is Fréchet differentiable on* $l_p$ *such that for* $\bar{x} = (\bar{t}_1, \bar{t}_2, \dots) \in l_p$, *the Fréchet derivative* $\nabla(P^m)(\bar{x}): l_p \to l_p$ *has the following representation.*

$$\nabla(P^m)(\bar{x}) = \begin{pmatrix} m\bar{t}_1^{m-1} & 0 & \dots \\ 0 & m\bar{t}_2^{m-1} & 0 \\ \vdots & 0 & \ddots \end{pmatrix}. \tag{8.16}$$

(ii) $P^m$ *is Gâteaux differentiable on* $l_p$ *such that for* $\bar{x} = (\bar{t}_1, \bar{t}_2, \dots) \in l_p$, *the Gâteaux derivative* $T'(\bar{x}): l_p \to l_p$ *satisfies that*

$$T'(\bar{x})(v) = \nabla(Q^m)(\bar{x})(v) = \{m\bar{t}_n^{m-1}t_n\}_{n=1}^{\infty}, \text{ for } v = \{t_n\}_{n=1}^{\infty} \in l_p \backslash \{\theta\}, \tag{8.17}$$

Now, let *m* be an odd number greater than 1 and we define a single-valued mapping $G: l_p \times l_p \to l_p$ as follows.

$$G(x, y) = P^m(x) + P^{m+1}(y), \text{ for any } x, y \in l_p. \tag{8.18}$$

**Example 8.6**. Let $1 < p < \infty$. Let $G: l_p \times l_p \to l_q$ be defined by (8.18). For any $\bar{x} = (\bar{t}_1, \bar{t}_2, \dots) \in l_p$ and $\bar{y} = (\bar{s}_1, \bar{s}_2, \dots) \in l_p$, by (8.15) with *m* being an odd number greater than 1, one sees that

(i) $G$ is order increasing with respect to the first variable at any given point $\bar{y} \in l_p$. By (8.16) or (8.17), we have

$$\nabla_1 G(\bar{x}, \bar{y}) = \nabla(P^m)(\bar{x}) = m \begin{pmatrix} \bar{t}_1^{m-1} & 0 & \dots \\ 0 & \bar{t}_2^{m-1} & 0 \\ \vdots & 0 & \ddots \end{pmatrix},$$

This immediately implies that

$$\nabla_1 G(\bar{x}, \bar{y})(v) \succcurlyeq_K \theta, \text{ for any } v \in K \backslash \{\theta\}.$$

(ii) $G$ is neither order increasing, nor order decreasing with respect to the second variable at any given $\bar{x} \in l_p$. By (8.16) or (8.17), we have

$$\nabla_2 G(\bar{x}, \bar{y}) = \nabla(P^{m+1})(\bar{y}) = (m+1) \begin{pmatrix} \bar{s}_1^m & 0 & \dots \\ 0 & \bar{s}_2^m & 0 \\ \vdots & 0 & \ddots \end{pmatrix}.$$

This immediately implies that, if $\bar{y} \neq \theta$, then

$$\nabla_2 F(\bar{x}, \bar{y})(u) \not\succcurlyeq_K \theta, \text{ for some } u \in K \backslash \{\theta\},$$

and

$$\nabla_2 F(\bar{x}, \bar{y})(w) \not\preccurlyeq_K \theta, \text{ for some } w \in K \backslash \{\theta\}.$$

□